\documentclass{amsart}
\usepackage{graphicx}
\usepackage{latexsym,amscd,amssymb,epsfig}

\theoremstyle{plain}
   \newtheorem{theorem}{Theorem}[section]
   \newtheorem{proposition}[theorem]{Proposition}
   \newtheorem{lemma}[theorem]{Lemma}
   \newtheorem{corollary}[theorem]{Corollary}
   
\theoremstyle{definition}
   \newtheorem{definition}[theorem]{Definition}
   \newtheorem{example}[theorem]{Example}
   
   \newtheorem{question}[theorem]{Question}
   \newtheorem{remark}[theorem]{Remark}

\numberwithin{equation}{section}

\newcommand\sym{{\mathfrak{S}}}
\newcommand\sh{{\mathrm{sh}}}
\newcommand\shf{{\mathrm{shf}}}
\newcommand\st{{\mathrm{std}}}

\newcommand\Knuth{{\underset{K}{\sim}}}
\newcommand\nuth{{\underset{K^*}{\sim}}}
\newcommand\Inv{{\mathrm{Inv}}}
\newcommand\Des{{\mathrm{Des}}}
\newcommand\Par{{\mathrm{Par}}}

\begin{document}
\title[Four orders on tableaux] {properties of four partial orders on standard
Young tableaux}

\author{ M\"{U}GE  TA\c{S}K{I}N }

\address{School of Mathematics\\
         University of Minnesota\\
         Minneapolis MN 55455}

\thanks{This research forms part of the author's doctoral thesis at the Univ. of
Minnesota, under the supervision of Victor Reiner, and partially supported by NSF grant
DMS-9877047.}

\begin{abstract}

Let $SYT_{n}$ be the set of all standard Young tableaux with $n$
cells. After recalling the definitions of four partial orders, the
weak, $KL$, geometric and chain orders on $SYT_n$ and  some of their
crucial properties, we prove three main results:

\begin{enumerate}
\item[$\bullet$]Intervals in any of these four orders essentially
describe the product in a Hopf algebra of tableaux defined by
Poirier and Reutenauer.

\item[$\bullet$] The map sending a tableau to its descent set induces
a homotopy equivalence of the proper parts of all of these orders on
tableaux with that of the Boolean algebra $2^{[n-1]}$. In
particular, the M\"obius function of these orders on tableaux is
$(-1)^{n-3}$.

\item[$\bullet$]  For two of the four orders, one can define a
more general order on skew tableaux having fixed inner boundary, and
similarly analyze their homotopy type and M\"obius function.
\end{enumerate}

\end{abstract}

\maketitle

\tableofcontents

\section{Introduction}

\begin{figure}
 \includegraphics[scale=0.50]{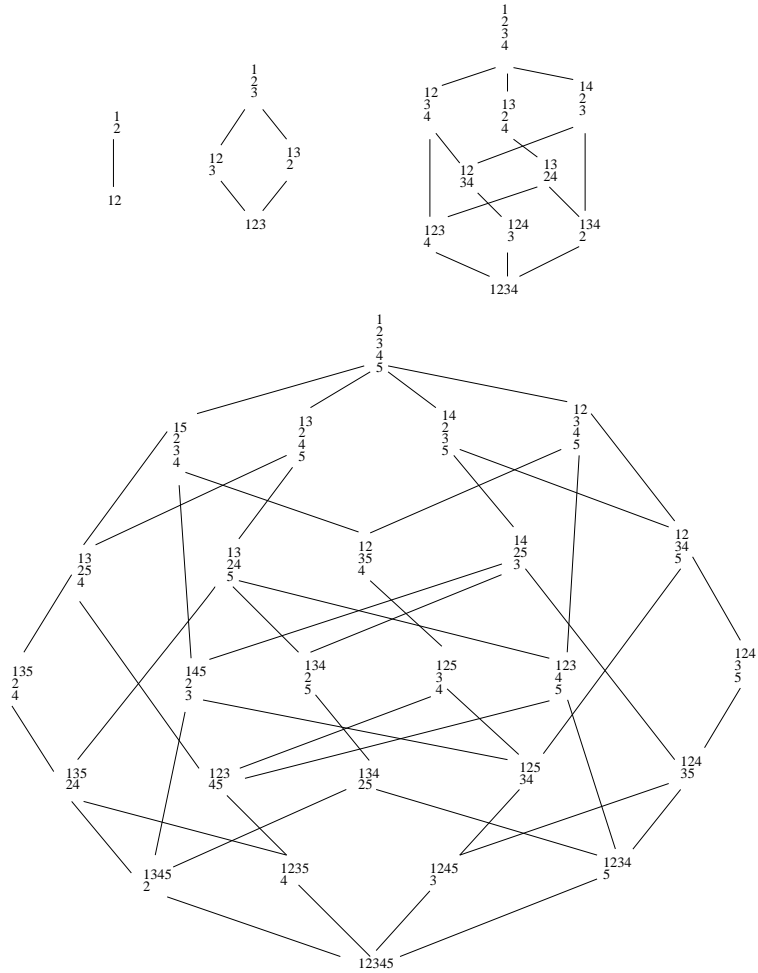} \caption{ \label{figure1}
 Chain, the weak,  $KL$
and geometric order  on $SYT_{n}$, which coincide for $n = 2,3,4,5$
(but not in general).}
\end{figure}

This paper is about four partial orders on the set $SYT_n$ of all
standard Young tableaux of size $n$ satisfying:
$$
\text{ weak order } \subsetneq \text{ Kazhdan-Lusztig (KL) order }
\subseteq \text{ geometric order } \subsetneq \text{ chain order }.
$$
Here $P \subset Q$ means that $u \leq v$ in $P$ implies $u \leq v$
in $Q$, in which case  we say $Q$ is {\it stronger} than $P$ (or $P$
is {\it weaker} then $Q$).

All four of these orders have appeared in the work of Melnikov
\cite{Melnikov1, Melnikov2,Melnikov3}, who refers to what we are
calling the weak order as the {\it induced Duflo order}. Roughly
speaking,

\begin{enumerate}
\item[$\bullet$]
the weak order is induced from the weak Bruhat order on the symmetric
group $\sym_n$ via the Robinson-Schensted insertion map,
\item[$\bullet$]
the KL order is induced by the Kazhdan-Lusztig preorder on $\sym_n$
arising in the theory of Kazhdan-Lusztig (right) cells,
\item[$\bullet$]
geometric order describes inclusions of certain algebraic varieties
indexed by tableaux ({\it orbital varieties}), and
\item[$\bullet$]
chain order is induced by the dominance order on partitions; for
each interval of values $[i,j]$, one restricts the tableau to these
values and compares the insertion shapes in dominance order.
\end{enumerate}

All four of these orders on $SYT_n$ coincide for $n \leq 5$, and are
depicted in Figure~\ref{figure1}.  For $n \geq 6$, they differ (see
Example~\ref{chain-order-shape-agree}  and
Example~\ref{KL-stronger-example}). After reviewing their
definitions in Section~\ref{Definitions-section}, we recall some of
their known properties in Section~\ref{known-properties-section}.

We then prove three main new results. The first result, proven in
Section~\ref{P-R-theorem-section}, relates to a Hopf algebra defined
by Poirier and Reutenauer \cite{Poirier-Reutenuer} whose basis
elements are indexed by standard Young tableaux $T$ of all sizes.
The multiplication in this Hopf algebra is somewhat nontrivial to
describe, but turns out to be described essentially by any of our
four partial orders\footnote{This result for the weak order was
asserted without proof in \cite[middle of p. 579]{Hivert-
Novelli-Thibon}.}.

\begin{theorem}
\label{P-R-theorem} For any of the four  partial order $\leq$ above,
one has
 $$
  T \ast S = \sum_{{\substack{R \in SYT_n:\\T/S \leq R \leq T\backslash S}}} R
 $$
where $T/S$ and $T\backslash S$ are obtained by sliding $S$ over $T$
from the left and from the bottom respectively.
\end{theorem}

\begin{figure}
\includegraphics[scale=0.45]{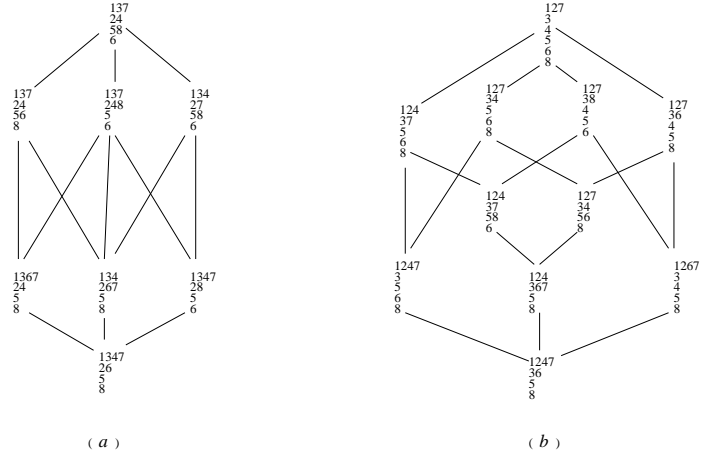}
\caption{\label{mob} ($a$) An interval in $(SYT_{8},\leq_{chain})$
and ($b$) an interval in $SYT_{8}$ ordered with $\leq_{weak}$,
$\leq_{KL}^{op}$ and $\leq_{geom}$ having M\"obius function $2$ and
$-2$ respectively.}
\end{figure}

The second result is about the M\"obius function and homotopy type
of these orders.  The weak Bruhat order on $\sym_n$ is well-known to
have each interval homotopy equivalent to either a sphere or a
point, and hence have M\"obius function values all in $\{\pm 1,0\}$.
Although it is not true in general for the intervals in  the weak,
$KL$, geometric and chain orders on $SYT_n$ (see Figure~\ref{mob}
for some examples)  the interval from bottom to top is homotopy
equivalent to either a sphere or a point. This result is  proven in
Section~\ref{homotopy-theorem-section}, by associating descent sets
to tableaux and thereby obtaining a poset map to a Boolean algebra.

\begin{theorem}
\label{homotopy-theorem} Let $\leq$ be any of the four partial
orders.
 Then the map $SYT_{n} \mapsto 2^{[n-1]}$ sending a
tableau to its descent set is order-preserving, and induces a
homotopy equivalence of the proper parts.

In particular, for any such order $\mu(\hat{0},\hat{1})=(-1)^{n-3}$.
\end{theorem}

The third result, proven in
Section~\ref{skew-homotopy-theorem-section} deals with a
generalization of the above orders to skew tableaux with fixed inner
boundary. The most crucial step in the proof is the application of
Rambau's Suspension Lemma \cite{Rambau} which makes the proof
(compared to  the standard methods in topological combinatorics)
much shorter and comprehensible. Given a partition $\mu$, let
$SYT^\mu_n$ denote the set of all skew standard tableaux of having
$n$ cells which are ``skewed by $\mu$'', that is, whose shape is
$\lambda / \mu$ for some $\lambda$. It turns out that two of the
four orders ($KL$, geometric) have a property (the {\it inner
translation} property; see Theorem~\ref{Vogan-involution}) which
allows us to generalize them on $SYT^\mu_n$.  Each of these skew
orders has a top element $\hat{1}$ and a bottom element $\hat{0}$,
so that one can speak of the homotopy of their {\it proper parts}
obtained by removing $\hat{0},\hat{1}$.

\begin{theorem}
\label{skew-homotopy-theorem} Let $\leq$ be $KL$ or geometric orders
on $SYT_n$.
 Then the associated order
$\leq$ on $SYT^\mu_n$ has the homotopy type of its proper part equal
to that of
$$
\begin{cases}
\text{ an } (n-2)-\text{dimensional sphere} & \text{ if }\mu \text{
is
rectangular,} \\
\text{ a point }  & \text{ otherwise.} \\
\end{cases}
$$
In particular, for any such order either
$\mu(\hat{0},\hat{1})=(-1)^{n-2}$ or $\mu(\hat{0},\hat{1})=0$,
depending on $\mu$.
\end{theorem}

\begin{center}
\begin{figure}
 \includegraphics[scale=0.60]{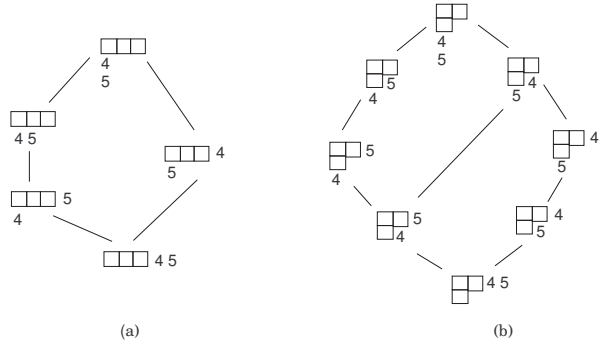}
 \vspace{-.3in}
 \caption{\label{skew-diagrams} An illustration of the skew orders on
$SYT_n^{\mu}$ for $n=2$. When  $\mu$  is ($a$) rectangular and ($b$)
nonrectangular, $SYT_2^{\mu}$ has its proper part homotopy
equivalent to a $0$-dimensional sphere and a  point respectively.}
\end{figure}
\end{center}

Figure~\ref{skew-diagrams} provides an illustration  for
Theorem~\ref{skew-homotopy-theorem}  where both posets are
considered  with $KL$ or geometric orders. In fact, this theorem
follows from a more general statement
(Proposition~\ref{skew-mobius-value-theorem}) about the homotopy
types of certain intervals, which applies to any order between the
weak and chain orders (including the weak order itself).

\vskip .1in We close this section with some context and motivation
for Theorems~\ref{P-R-theorem} and \ref{homotopy-theorem}, stemming
from two commutative diagrams that appear in the work of Loday and
Ronco \cite{Loday-Ronco}
\begin{equation}
\label{triangles}
\begin{array}{cll}
        \sym_{n} & \longrightarrow & Y_{n}\\
     &  \searrow   & \downarrow  \\
                           &       &    2^{[n-1]} \\
\end{array} \hskip .4in
\begin{array}{cll}
        \mathbb{Z}\sym & \longleftarrow & \mathbb{Z}Y \\
     &  \nwarrow   & \uparrow  \\
                           &       &    \Sigma  \\
\end{array}
\end{equation}

In the left diagram of \eqref{triangles}, $Y_n$ denotes the set of
planar binary trees with $n$ vertices. The horizontal map sends a
permutation $w$ to a certain tree $T(w)$, and has been considered in
many contexts (see e.g. \cite[\S 1.3]{Stanley},
\cite[\S9]{Bjorner4}).  The southeast map $\sym_n \rightarrow
2^{[n-1]}$ sends a permutation $w$ to its descent set $\Des_L(w)$.
These maps of sets become order-preserving if one orders $\sym_n$ by
weak order, $Y_n$ by the {\it Tamari order} (see
\cite[\S9]{Bjorner4}), and $2^{[n-1]}$ by inclusion.

Indeed, the  order preserving maps of the first diagram induce the
inclusions of Hopf algebras in the second diagram of
\eqref{triangles},  in which $\mathbb{Z}\sym$ is the
Malvenuto-Reutenauer algebra, $\mathbb{Z}Y$ is a subalgebra
isomorphic to Loday and Ronco's free {\it dendriform algebra} on one
generator \cite{Loday-Ronco1}, and $\Sigma$ is a subalgebra known as
the algebra of {\it noncommutative symmetric functions}. In
\cite{Loday-Ronco}, Loday and Ronco proved a description of the
product structure for each of these three algebras very much
analogous to Theorem~\ref{P-R-theorem}, which should be viewed as
the analogue replacing $\mathbb{Z}Y$ by $\mathbb{Z}SYT$; see Theorem
\ref{Loday-Ronco-Malvenuto-Reutenauer} below  for their description
of the product in $\mathbb{Z}\sym$. The analogy between the standard
Young tableaux $SYT_n$ and the planar binary trees $Y_n$ is
tightened further by recent work of Hivert, Novelli and Thibon
\cite{ Hivert- Novelli-Thibon}.  They show that the planar binary
trees $Y_n$ can be interpreted as the plactic monoid structure given
by a Knuth-like relation, similar to the interpretation of the set
of standard Young tableaux as Knuth/plactic classes.

We were further motivated in proving Theorem~\ref{P-R-theorem} by
the results of Aguiar and Sottile in  \cite{Aguiar-Sottile} and
\cite{Aguiar-Sottile2} where the M\"obius functions of the weak
order on $\sym_n$ and Tamari order on $Y_n$ have  key roles in
understanding the structures of the  Hopf algebras of permutations
and planar binary trees.

In \cite[Remark 9.12]{Bjorner4}, Bj\"orner and Wachs (essentially)
show that the triangle on the left induces a diagram of homotopy
equivalences on the proper parts of the posets involved.
Theorem~\ref{homotopy-theorem} gives the analogue of this statement
in which one replaces $(Y_n, \leq_{Tamari})$ by $(SYT_n,\leq)$ where
$\leq$ is any order between the weak and chain orders.

\vskip .3in
\section{Definitions}
\label{Definitions-section}

\vskip .1in
\subsection{Chain order}
\label{chain-order-subsection} \noindent

 The first partial order on $SYT_n$ that will be discussed  is the
 strongest one:  {\it chain order}.

  Given  $T \in SYT_n$, we  denote by $\sh(T)$ the partition corresponding
to the shape of $T$. For $1\leq i< j \leq n$, let  $T_{[i,j]}$ be
the skew subtableau obtained by restricting $T$ to the segment
$[i,j]$. Let $\st(T_{[i,j]})$ be  the tableau obtained by  lowering
all entries of $T_{[i,j]}$  by $i-1$ and sliding it into normal
shape by jeu-de-taquin \cite{Schutzenberger2}.

    The definition of chain order also involves the dominance order.
We denote by $(\Par_n, \leq^{op}_{dom})$ the set of all partitions
of the number $n$ ordered by the {\it opposite (or dual) dominance
order}, that is, $\lambda \leq^{op}_{dom} \mu$ if
$$
\lambda_1 + \cdots + \lambda_k \geq \mu_1 + \cdots + \mu_k \text{
for all } k.
$$

\begin{definition}
\label{chain-order-def} Let $S, T \in SYT_n$ and  We say $S$ is less
that $T$ in chain order ($S\leq_{chain} T$) if
 for  every $1\leq i< j \leq n$,
 $$\sh(\st(S_{[i,j]}))~\leq^{op}_{dom}~\sh(\st(T_{[i,j]})).$$
\end{definition}

\vskip .2in
\subsection{Weak order}
\label{weak-order-subsection}
 \noindent

 Before giving the definition of the weak order   it is
necessary to recall the Robinson-Schensted  $(RSK)$ correspondence;
see \cite[\S 3]{Sagan} for more details and references on $RSK$.
The $RSK$ correspondence is a bijection between $\sym_{n}$ and
 $\{(P,Q): P, Q \in  SYT_{n} \mbox{ of same shape}\}$.
Here $P$ and $Q$ are called  the {\it insertion} and {\it recording
tableau} respectively. Knuth \cite{Knuth1} defined an equivalence
relation $\Knuth$ on $\sym_n$ with the property that ~$u \Knuth w$~
if and only if they have the same insertion tableaux $P(u)=P(w)$. We
will denote the corresponding equivalence classes in $\sym_n$ by
$\{C_T\}_{T\in SYT_n}$.

We now recall the ({\it right}) {\it weak} {\it Bruhat order},
$\leq_{weak}$, on  $\sym_n$.  It is the transitive closure of the
relation $u \leq_{weak} w$ if $w=u \cdot s_i$ for some $i$ with
$u_{i} < u_{i+1}$, and where $s_i$ is the adjacent transposition $(i
\,\, i+1)$.  The weak order has an alternative characterization
\cite[Prop. 3.1]{Bjorner2}  in terms of {\it (left) inversion sets}
$$
\Inv_L(u):=\{(i,j): 1 \leq i < j \leq n ~\text{ and }~
         u^{-1}(i) > u^{-1}(j)\},
$$
namely ~$u \leq_{weak} w$~ if and only if $\Inv_L(u)\subset
\Inv_L(w)$.

\vskip .1in
For $1\leq i< j \leq n$ let $[i,j]$ be a
  segment of the alphabet $[n]$ and  $u_{[i,j]}$ be the
subword of $u$ obtained by restricting to the alphabets $[i,j]$ and
$\st(u_{[i,j]})$ in $\sym_{j-i+1}$ be the word obtained from
$u_{[i,j]}$ by subtracting $i-1$ from each letter.

 In fact $\Inv_L(u)\subset \Inv_L(w)$ gives
$\Inv_L(u_{[i,j]})\subset \Inv_L(w_{[i,j]})$ for all $1\leq i < j
\leq n$ and  hence
\begin{equation}\label{welldefined-weak-order1}
u \leq_{weak} w ~\mbox{ implies } ~u_{[i,j]} \leq_{weak} w_{[i,j]}
~\mbox{ for all } ~1\leq i<j\leq n.
\end{equation}

The following   basic fact about $RSK$, Knuth equivalence, and
jeu-de-taquin  are essentially due to Knuth and Sch\"utzenberger;
see Knuth \cite[Section 5.1.4]{Knuth} for detailed explanations.

\begin{lemma}
\label{j-d-t-initial-final} Given $u \in \sym_n$,
 let $P(u)$ be the insertion tableau of $u$.  Then for $1\leq i<j \leq
 n$,
$$
\st(P(u)_{[i,j]}) = P(\st(u_{[i,j]})).$$
\end{lemma}
\vskip .1in

Furthermore one can use  Greene's theorem \cite{Greene} for
the following fact:
\begin{equation}\label{welldefined-weak-order2}
\mbox{If }u \leq_{weak} w ~\mbox{ then }
~\sh(\st(P(u)_{[i,j]}))\leq^{op}_{dom}\sh(\st(P(u)_{[i,j]})) ~\mbox{
for all } ~1\leq i<j\leq n.
\end{equation}
Now  \eqref{welldefined-weak-order1} and
\eqref{welldefined-weak-order2} shows that the following order is
weaker than  chain order on $SYT_n$ and hence it is well defined.

\begin{definition}
\label{weak-order-def}
 The {\it weak order} $(SYT_n,\leq_{weak})$,
first introduced by Melnikov \cite{Melnikov1} under the name {\it
induced Duflo order}, is the partial order induced by taking
transitive closure of  the following rule. Denoting  the Knuth class
of $T$ by $C_T$,
$$
\begin{aligned}
S\leq_{weak} T & \mbox{ if there exist } ~\sigma\in C_S, ~\tau\in C_T \\
               & \mbox{ such that } ~\sigma\leq_{weak} \tau.
\end{aligned}
$$

\end{definition}

 The necessity of taking the transitive closure in the definition of the
 weak order is illustrated by the following example (cf.
Melnikov \cite[Example 4.3.1]{Melnikov1}).

\begin{example}
Let
$R=\scriptstyle{\begin{array}{ccc} 1& 2 & 5 \\  3& 4 &
\end{array}}$,
\hskip .05in $S=\scriptstyle{\begin{array}{ccc} 1& 4 & 5 \\  2&  &
\\ 3 & & \end{array}}$,
\hskip .05in $T=\scriptstyle{
\begin{array}{ccc} 1 & 4 &  \\  2 & 5 &  \\ 3 & &\end{array}}$
with
$$
 \begin{array}{ll}
       &  C_{R}=\{ 31425, 34125, 31452, 34152, 34512 \}, \\
       & C_{S}=\{ 32145, 32415, 32451, 34215, 34251 ,34521 \}, \\
       & C_{T}= \{ 32154, 32514, 35214, 32541, 35241\}. \\
 \end{array}
 $$
Here $R <_{weak} S $ since $34125 <_{weak} 34215 =34125 \cdot
s_{3}$, and $S <_{weak} T $ since $32145 <_{weak} 32154 = 32145
\cdot s_{4}$. Therefore $ R <_{weak} T $.

    On the other hand, for every $
\rho \in C_{R}$ one has $(2,4) \in \Inv_L(\rho)$, whereas for every
$\tau \in C_{T}$ one has $(2,4) \notin \Inv_{L}(\tau)$. This shows
that there is no $\rho \in C_{R}$ and $\tau \in C_{T}$ such that
$\rho <_{weak} \tau$.
\end{example}

\vskip .2in
\subsection{Kazhdan-Lusztig order}
\label{KL-order-subsection} \noindent

It turns out that $RSK$ is closely related to \emph{Kazhdan-Lusztig}
preorders  on $\sym_n$. Recall that a {\it preorder} on a set $X$ is
a binary relation $\leq$ which is reflexive ($x \leq x$) and
transitive ($x \leq y, y \leq z$ implies $x \leq z$).  It need not
be antisymmetric, that is, the equivalence relation $x \sim y$
defined by $x \leq y, y \leq x$ need not have singleton equivalence
classes. Note that a preorder induces a {\it partial order} on the
set $X/_{\sim}$ of equivalence classes.

Kazhdan and Lusztig \cite{Kazhdan-Lusztig} introduced two preorders
(the left and right $KL$ preorders) on Coxeter groups whose
equivalence classes are called the left and  right cells
repectively. The theory of left (or right) cells provides a
decomposition of the regular representation of the Hecke algebras of
Coxeter groups (c.f. \cite[Chapter 6]{Bjorner-Brenti}) such that, in
case the Coxeter group is  $\sym_n$, each summand is irreducible.

 In this paper we will denote by $\leq^{op}_{KL}$
the {\it opposite} of the usual $KL$ right preorder on $\sym_n$. For
example, with our convention, the identity element $1$ and the
longest element $w_0$ satisfy $1 \leq^{op}_{KL} w_0$. It turns out
\cite{Kazhdan-Lusztig} (and explicitly in \cite[p.
54]{Garsia-McLarnan}) that the associated equivalence relation for
this $KL$ preorder is the Knuth equivalence~$\Knuth$. Hence an
equivalence class (usually called either a {\it Knuth class} or {\it
plactic class} or a {\it Kazhdan-Lusztig right cell} in $\sym_n$)
corresponds to a tableau $T$ in $SYT_n$.

\begin{definition}
\label{KL-order-def}
 {\it $KL$ order} on $SYT_n$   is defined by the rule
$$
 S\leq^{op}_{KL} T \hskip .1in \text{if } \hskip .1in C_S \leq^{op}_{KL} C_T
$$
 where $C_S$ is the Knuth class (or
$KL$ right cell) in $\sym_n$ corresponding to  $S \in SYT_n$.
\end{definition}

For later use, we now  recall the basic construction of the $KL$
right preorder on $\sym_n$. Recall that the {\it right descent set}
$D_R(u)$ and  the {\it left descent set} $D_L(u)$ of a permutation
$u \in \sym_n$,  are defined by
$$
\begin{aligned}
\Des_{R}(u) & :=
   \{(i,i+1): 1 \leq i \leq n-1 ~\text{ and }~
       u(i) > u(i+1)  \} \\
\Des_{L}(u) & :=
   \{(i,i+1): 1 \leq i \leq n-1 ~\text{ and }~
       u^{-1}(i) > u^{-1}(i+1)  \} \\
            & = \Inv_L(u) \cap S\\
\end{aligned}
$$
 where $S=\{(i,i+1): 1 \leq i \leq n-1 \}$.  In what follows, we
will often identify the set $S$ of adjacent transpositions with the
numbers $[n-1]:=\{1,2,\ldots,n-1\}$ via the obvious map $(i,i+1)
\mapsto i$.

\vskip .1in In \cite{Kazhdan-Lusztig}, Kazhdan and Lusztig prove the
existence of unique polynomials $\{P_{u,w}(q)\} \subseteq
\mathbb{Z}[q]$ indexed by permutations in $\sym_n$. Denoting by
~$\leq$~  the Bruhat order on $\sym_n$, $l(u)$ the length of the
permutation  $u$  and $l(u,w)=l(w)-l(u)$, these polynomials satisfy:
\begin{equation}\label{basic-properties-KL-polynomial}
\begin{aligned}
&P_{u,w}(q)=1~\text{ if }~u=w,\\
& P_{u,w}(q)=0 ~\text{ if }~u \nleq w,  \\
& \mathrm{deg}(P_{u,w}(q)) \leq \frac{1}{2}(l(u,w)-1).
\end{aligned}
\end{equation}

Let $[q^k]P_{u,w}(q)$ denote  the coefficient of $q^k$ in
$P_{u,w}(q)$ and   define
\begin{equation}\label{KL-Mu-function}
\overline{\mu}(u,w):=
\begin{cases}
[q^{\frac{l(u,w)-1}{2}}]P_{u,sw}(q)&~~ \text{ if } l(u,w)
\text{ is odd }\\
~~0&~~  \text{ otherwise. }
\end{cases}
\end{equation}
Then a   recursive formula for these polynomials is given in the
following way: For $u \leq w$ and $s \in D_L(w)$,

\begin{equation}\label{recursive-formula-KL-polynomial}
P_{u,w}(q)= q^{1-c} P_{su,sw}(q)~+~q^c P_{u,sw}(q)~-~ \sum_{\{v:s\in
D_L(v)\}} q^{l(v,w)/2}~ \overline{\mu}(v,sw)~P_{u,v}(q)
\end{equation}
where $c=1$ if  $s\in D_L(u)$ and $c=0$ otherwise. Moreover the dual
of right $KL$ preorder on $\sym_n$ is given by taking the transitive
closure of the following  relation:
\begin{equation}\label{KL-preorder}
u \leq^{op}_{KL} w~~~ \text{ if }~~~
\begin{cases}
D_R(w)- D_R(u) \not = \varnothing,\\
 \text{ and }  \\
 \overline{\mu}(u,w)\not=0~
~\text{ or }~ \overline{\mu}(w,u)\not=0.
\end{cases}
\end{equation}

\vskip .2in
\subsection{Geometric order}
\label{geometric-order-subsection}\noindent

The final order on $SYT_n$ to be discussed in this paper relates to
the preorder on $\sym_n$ induced from geometric order on the orbital
varieties associated to the Lie algebra $\mathfrak{sl}_n$.

The theory of  orbital varieties arise from the work of N.
Spaltenstein \cite{Spaltenstein1,Spaltenstein2} and R. Steinberg
\cite{Steinberg1,Steinberg2} on the unipotent variety of a connected
complex semi-simple group $G$. They have a key role  in the studies
of primitive ideals (i.e. annihilators of irreducible
representations) in the enveloping algebra
$\mathfrak{U}(\mathfrak{g})$ of Lie algebra $\mathfrak{g}$
corresponding to $G$ (c.f. \cite{Melnikov4},
\cite{Collingwood-MacGovern}, \cite{Leeuwen2}). They also play an
important role in Springer's Weyl group representations.

Let $\mathfrak{g}$ be the Lie algebra of $G$ and  $B$ be the Borel
subgroup of $G$ given with respect to some triangular decomposition
$\mathfrak{g}= \mathfrak{n}^-  \oplus \mathfrak{h} \oplus
\mathfrak{n}$ such that $\mathfrak{h}$ is a Cartan subalgebra and
$\mathfrak{n}$ is the corresponding nilradical.

For given $\eta \in \mathfrak{n}$, we denote by $\mathcal{O}_\eta$
the nilpotent orbit determined  by  the adjoint action of $G$ on
$\eta$. Therefore $\overline{\mathcal{O}}_\eta$ is an irreducible
variety. Now an {\it orbital variety} $\mathcal{V}$  associated to
$\mathcal{O}_\eta$ is defined to be an irreducible component of the
intersection $\mathcal{O}_{\eta} \cap \mathfrak{n}$. Given orbital
varieties $\mathcal{V}$ and $\mathcal{W}$, the {\it geometric order}
is defined by
$$
\mathcal{V}\leq_{geom}\mathcal{W} ~~\text{ if }~~\mathcal{W}
\subseteq \overline{\mathcal{V}}
$$
 where $\overline{\mathcal{V}}$ denotes the
Zariski closure of $\mathcal{V}$ inside $\mathfrak{n}$. The only
general description of orbital varieties provided below is due
Steinberg \cite{Steinberg1}.

\vskip .1in
  Given a positive root system  $R^+\subset\mathfrak{h}^*$,
 recall that $\mathfrak{n}= \oplus_{\alpha \in
R^+} \mathfrak{g}_\alpha $ where  $\mathfrak{g}_\alpha$ is the root
space corresponding to $\alpha$. Let $W$ be the Weyl group of
$\mathfrak{g}$ generated by simple roots in $R^+$, and  for  $w \in
W$ let
$$
\mathfrak{n}\cap^{w} \mathfrak{n}:= \bigoplus_{\alpha \in ~w(R^+)
\cap R^+} \mathfrak{g}_{\alpha}.
$$

Since  $B$ is an irreducible closed subgroup   of $G$, the action of
$B$ on ~$\mathfrak{n}\cap^{w} \mathfrak{n}$~ gives an irreducible
locally closed subvariety $B(\mathfrak{n}\cap^{w} \mathfrak{n})$
which, therefore, lies in a unique nilpotent orbit
$\mathcal{O}_\eta$  for some $\eta \in \mathfrak{n}$ and
$\overline{G(\mathfrak{n}\cap^{w}
\mathfrak{n})}=\overline{\mathcal{O}}_{\eta}$.
 By the result of Steinberg
\begin{equation}
\label{orbital-variety}
 \mathcal{V}_{w}:=
\overline{B(\mathfrak{n}\cap^{w}\mathfrak{n})}\cap
\mathcal{O}_{\eta}
\end{equation}
is an orbital variety and the map $w \mapsto \mathcal{V}_w$ is a
surjection. Moreover geometric order induces a preorder on  $W$ such
that, for $u,w \in W$
\begin{equation}
\label{geometric-order-on-permutations}
\begin{aligned}
u \leq_{geom} w ~\text{ if }~& \mathcal{V}_w
\subset\overline{\mathcal{V}}_u \\
& ~\text{ or equivalently }~
\overline{B(\mathfrak{n}\cap^{w}\mathfrak{n})}
\subseteq\overline{B(\mathfrak{n}\
\cap^{u}\mathfrak{n})}.
\end{aligned}
\end{equation}
\vskip .1in
 According to Steinberg \cite{Steinberg1}, the fibers of the map $w \mapsto
\mathcal{V}_w$ for
 $\mathfrak{g}=\mathfrak{sl}_n$ are the Knuth classes of $\sym_n$  and
therefore each orbital variety $\mathcal{V}$ in $\mathfrak{sl}_n$
can be identified with some $T \in SYT_n$ i.e.,
$\mathcal{V}=\mathcal{V}_T$.
 This leads to  the following definition.

\begin{definition}
The  {\it geometric order} on $SYT_n$, $(SYT_n,\leq_{geom})$, is
given by the following rule:
$$
 S\leq_{geom} T \hskip .1in \text{ if } \hskip .1in  \mathcal{V}_T \subset
\overline{\mathcal{V}}_S.
 $$
\end{definition}

When $\mathfrak{g}=\mathfrak{sl}_n$, an explicit description of
orbital varieties can be given in the following way. Let $B$ to be
the Borel subgroup of invertible upper triangular $n\times n$
matrices given by   the Cartan decomposition of $\mathfrak{g}$ with
 Cartan subalgebra $\mathfrak{h}$ of trace $0$ diagonal matrices
and nilradicals $\mathfrak{n}$ and $\mathfrak{n}^-$,  whose elements
are strictly upper and  strictly lower triangular matrices
respectively. Then the set of matrices $\{E_{ij}\}_{i<j}$
(and~$\{E_{ij}\}_{i>j}$), where $E_{i,j}$ has $1$ on the  position
$(i,j)$ and $0$ elsewhere, provides a basis
 for $\mathfrak{n}$ (respectively~$\mathfrak{n^-}$).

The  action of the Weyl group $\sym_n$ on $E_{i,j}$ can be described
by
     $$
     w\cdot E_{i,j}= p_{w}~ E_{i,j}~
     p_{w}^{-1}=E_{w(i),w(j)}
     $$
where $p_{w}$ is the permutation matrix of $w \in \sym_n$ and this
leads to the following characterization
\begin{equation}
\label{geom-eq1}
\mathfrak{n}\cap^{w}\mathfrak{n}=\mathrm{span}\{E_{i,j}\mid (i,j)
\not \in \Inv_L(w)\}.
\end{equation}
On the other hand  the adjoint action of $B$ on $E_{i,j}$ sweeps the
corner at $(i,j)$ to the northeast direction. In other words $B\cdot
E_{i,j}$ consists of  all  matrices of rank $1$, having  a nonzero
entry at $(i,j)$ and  all other  nonzero entries are located at some
positions to the northeast of $(i,j)$. Therefore all matrices in
$\overline{B(\mathfrak{n}\cap^{w}\mathfrak{n})}$ have their nonzero
entries in some boundary provided by $\{B \cdot E_{i,j}\mid(i,j)
\not \in \Inv_L(w)\}$, and
$\mathcal{V}_w=\overline{B(\mathfrak{n}\cap^{w}\mathfrak{n})}\cap
\mathcal{O}_{\eta}$ consists of all those matrices in
$\overline{B(\mathfrak{n}\cap^{w}\mathfrak{n})}$ whose Jordan form
is the same as  that of $\eta$. Recall that $\eta$ is uniquely
determined by the condition $\overline{G(\mathfrak{n}\cap^{w}
\mathfrak{n})}=\overline{\mathcal{O}}_{\eta}$.  Actually one can
show that the partition  determined  by the Jordan form of $\eta$
and the partition obtained from $w$ through the $RSK$ correspondence
are the same.

\vskip .1in
There is also a bijection, revealed by Steinberg
\cite{Steinberg2}, between the orbital varieties determined by
$\eta$ and Springer fiber  $\mathcal{F}_{\eta}$ of the complete flag
variety $\mathcal{F}$. Moreover geometric order results in an
ordering between the irreducible components of $\mathcal{F}_{\eta}$.
We next discuss this connection.

Let $\lambda=J(\eta)$ be the Jordan form of $\eta$,
$\mathcal{O}_\lambda=\{\eta \mid J(\eta)=\lambda\}$ be the
$GL(V)$-orbit of $\eta$ and
$$
\widetilde{\mathcal{O}}_{\lambda}:=\{(\eta,f) \mid \eta \in
\mathcal{O}_{\lambda}, f \in \mathcal{F}: \eta(f)\subset f \}.
$$
Here $GL(V)$ acts on $\mathcal{O}_{\lambda}$ and $\mathcal{F}$ by
conjugation and left translation respectively; therefore it acts on
$\widetilde{\mathcal{O}}_{\lambda}$,  and the projections onto
$\mathcal{O}_{\lambda}$ and  $\mathcal{F}$ are equivariant maps. We
have the following diagram:
$$
\begin{array}{lrclr}
       & &  \widetilde{\mathcal{O}}_{\lambda}&& \\
     &       \swarrow &  &\searrow&\\
                      \mathcal{O}_{\lambda}      &       &&& \mathcal{F}  \\
\end{array}
$$
In this diagram,  the fiber  of any $\eta \in \mathcal{O}_{\lambda}$
is equal to $\mathcal{F}_\eta:=\{(\eta,f): f \in \mathcal{F},
\eta(f)\subset f\}$. Since $GL(V)$ is irreducible and its action on
$ \mathcal{O}_\lambda$ is transitive, the irreducible components of
this {\it Springer fiber} $\mathcal{F}_{\eta}$  are in bijection
with the irreducible components of
$\widetilde{\mathcal{O}}_{\lambda}$. On the other hand for any $f
\in \mathcal{F}$, let $B$ be the Borel subgroup of $GL(V)$ which
fixes $f$ and let $\mathfrak{n}$ be nilradical of the corresponding
Borel algebra $\mathfrak{b}$. Then the fiber of $f$ is equal to
$\{(\eta,f): \eta \in \mathcal{O}_\lambda\cap \mathfrak{n}\}$ and
again the transitivity of  the action of $GL(V)$ on $\mathcal{F}$
implies that  the irreducible components of $\mathcal{O}_\lambda\cap
\mathfrak{n}$ are in bijection with the irreducible components of
$\widetilde{\mathcal{O}}_{\lambda}$. These two bijections determine
the correspondence between the orbital varieties and the irreducible
components of Springer fibers in the flag variety. The geometric
order describes the inclusions among (the closures of) these
components as one varies $\lambda$, in either context.

\vskip .3in
\section{Known properties}
\label{known-properties-section}
\vskip .1in
 In this section we recall  some of
the main properties of these four orders which we need later in
proving our main results. These properties also can be found in or
deduced from the works of Melnikov
\cite{Melnikov1,Melnikov2,Melnikov4} and Barbash and Vogan
\cite{Barbash-Vogan}. In order to make these posets more
understandable we  provide the proofs of those which are
combinatorially approachable, while for those which need theoretical
approaches the reader is directed to the references.

\vskip .1in
\subsection{Restriction to segments}
\label{Restriction-to-segments-section} \noindent

For $u \in \sym_n$ and $T \in SYT_n$ recall the definitions of
~~$\st(u_{[i,j]})$ and $\st(T_{[i,j]})$ from
Section~\ref{weak-order-subsection} and
Section~\ref{chain-order-subsection} respectively.
\vskip .1in Say
that a family of preorders $\leq$ on $\sym_n$  {\it restricts to
 segments} if
$$
u \leq w ~\text{  implies }~ ~\st(u_{[i,j]}) \leq
 \st(w_{[i,j]})~ \text{  for all } ~1\leq i<j \leq n.
$$

Melnikov shows in \cite[Page 45]{Melnikov1} the preorder
$\leq_{geom}$ on the Weyl group $W$  of any reductive Lie algebra
restricts to $W_I$, where  $I$ is any subset of simple roots
generating $W_I$. Therefore geometric order on $\sym_n$ restricts to
segments. The same fact about $KL$ preorder was first shown by
Barbash and Vogan \cite{Barbash-Vogan} for arbitrary finite Weyl
groups (see also work by Lusztig \cite{Lusztig2}) whereas the
generalization to Coxeter groups  is due to Geck \cite[Corollary
3.4]{Geck}. On the other hand this result for the weak order on
$\sym_n$ follows from an easy analysis on the (left) inversion sets.

\vskip .1in
We say the order  $\leq$ on $SYT_n$  {\it restricts to
segments} if
 $$S \leq T~~\text{ implies }~~\st(S_{[i,j]}) \leq
 \st(T_{[i,j]})~~ \text{ for all } ~~1\leq i<j \leq n.
 $$
The following result for  the weak, $KL$ and geometric order on
$SYT_n$  is an easy consequence of the above discussion together
with Lemma \ref{j-d-t-initial-final}, whereas for chain order it
follows directly from its definition.

\begin{corollary}
\label{segment-restriction-lemma} On $SYT_n$ all of the four orders
restrict to  segments of standard Young tableaux.
\end{corollary}

\vskip .1in
In fact  any order ~$\leq$~ on
$SYT_n$ which is stronger than the weak order and which restricts to
segments shares a crucial property that we describe now.

Recall that {\it (left) descent set} of a permutation $\tau$ is
defined by
$$
\Des_{L}(\tau) :=\{i: 1 \leq i \leq n-1 ~\text{ and }~
       \tau^{-1}(i) > \tau^{-1}(i+1)  \} \\
$$
 As a consequence of a
well-known properties of $RSK$, the left descent set $\Des_L(-)$ is
constant on Knuth classes $C_T$; the {\it descent set} of the
standard Young tableau $T$ is described intrinsically by
$$
\begin{aligned}
\Des(T) &:=
   \{(i,i+1): 1 \leq i \leq n-1 \text{ and }  \\
   & \qquad \qquad  i+1 \text{ appears in a row below } i \text{ in }T\}.
\end{aligned}
$$

We let $(2^{[n-1]}, \subseteq)$ be the Boolean algebra of all
subsets of $[n-1]$ ordered by inclusion. \vskip .2in

\begin{lemma}\label{tableau-descent-lemma}
Let $\leq$ be any order on $SYT_n$ which is stronger than the weak
order and restricts to segments. Then the map
$$(SYT_n, \leq) \mapsto (2^{[n-1]}, \subseteq)$$
sending any tableau $T$ to its descent set $\Des(T)$ is order
preserving.
\end{lemma}
\begin{proof} For $n=2$, such an order is isomorphic to weak order
on $SYT_2$ and the  statement follows directly by examination of
Figure~\ref{figure1}. For $n>2$, one can use the fact that
$$\Des_L(T)=\Des_L(T_{[1,n-1]})\cup\Des_L(T_{[2,n]})$$
to get the desired result by induction.
\end{proof}

\vskip .2in
\subsection{Poset morphisms}
\label{Poset-morphisms-section} \noindent

For the record, we note here some symmetries and order-preserving
maps of $\leq_{chain}$, $\leq_{weak}$,  $\leq^{op}_{KL}$ and
$\leq_{geom}$ on $SYT_n$,
 to other posets.

\begin{proposition}
\label{order-preserving-maps}
 Let $\leq$ represent to  any of the
orders  $\leq_{weak}$, $\leq^{op}_{KL}$
 $\leq_{geom}$ or $\leq_{chain}$
 on $SYT_n$. Then the following maps are order
preserving:
\begin{enumerate}
\item[(i)]  The map
$$
(SYT_{n},\leq) \rightarrow  (2^{[n-1]}, \subseteq)
$$
sending a tableau $T$ to its descent set $\Des(T)$.
\item[(ii)] The map
$$
(SYT_{n},\leq) \rightarrow (\Par_{n},\leq^{op}_{dom})
$$
sending $T$ to its shape $\lambda(T)$.
\end{enumerate}
On the other hand for $\leq$ equal to  any of the orders
$\leq_{weak}$, $\leq^{op}_{KL}$,  $\leq_{geom}$ or $\leq_{chain}$
\begin{enumerate}
\item[(iii)] the Sch\"utzenberger's evacuation map
$$
(SYT_{n},\leq) \rightarrow  (SYT_{n},\leq)
$$
sending $T$ to its evacuation tableau $T^{evac}$ is an poset
automorphism,
\end{enumerate}
 whereas for  $\leq$ equal to
 $\leq_{weak}$,  $\leq^{op}_{KL}$ or $\leq_{chain}$
\begin{enumerate}
\item[(iv)] the map
$$
(SYT_{n},\leq) \rightarrow  (SYT_{n},\leq)
$$
sending $T$ to its transpose $T^t$ is a poset anti-automorphism.
\end{enumerate}
\end{proposition}

\begin{proof}
The first assertion follows from Lemma~\ref{tableau-descent-lemma},
since all of the four orders are stronger than the weak order and
restrict to segments.

  Second assertion for $\leq_{chain}$ follows from its definition. For
$\leq_{weak}$, as it mentioned earlier,  one can apply Greene's
Theorem \cite{Greene}. If $S \leq_{geom} T$ then  there are orbital
varieties given by $\mathcal{V}_S$ and $\mathcal{V}_T$ such that
$\mathcal{V}_T\subseteq\overline{\mathcal{V}}_S$. Now the nilpotent
orbits that these orbital varieties  belong to can be characterized
by the partition given by $\sh(S)$ and $\sh(T)$. Moreover we have
$\mathcal{O}_{\sh(T)}\subseteq\overline{\mathcal{O}}_{\sh(S)}$. By
the result of Gerstenhaber, see \cite[Chapter
6]{Collingwood-MacGovern} for example, last inclusions implies
$\sh(T)\leq_{dom} \sh(S)$, proving the statement for geometric
order. For $KL$ order  the proof based on the theory that relates
the Kazhdan-Lusztig cells to the primitive ideals: let
$\mathfrak{g}$ be a semisimple algebra with universal enveloping Lie
algebra $U(\mathfrak{g})$ and Weyl group $W$. As it is  shown in
\cite{Barbash-Vogan} and \cite{Borho1}, for any primitive ideal $I$
of $U(\mathfrak{g})$,  the set of the form $\{w \in W \mid
I_{w}=I\}$ can be characterized by a Kazhdan-Lusztig left cell.
Moreover  $v <_{KL}^{op} w$ (right dual $KL$ order) if and only if
$I_{v^{-1}} \subseteq I_{w^{-1}}$, whence the associated variety of
the primitive ideal $I_{w^{-1}}$ is contained in that of
$I_{v^{-1}}$. On the other hand by the result of Borho and Brylinski
\cite{Borho1} and Joseph \cite{Joseph} associated variety of a
primitive ideal  is the closure of a nilpotent orbit in
$\mathfrak{g}^*$.  In our case  $\mathfrak{g}=\mathfrak{sl}_n$,
$W=\sym_n$ and the nilpotent orbits are characterized by partitions
of $n$, therefore  the result of Gerstenhaber reveals the desired
property on the shapes of the corresponding tableaux of $v$ and $w$.

  The assertions about transposition and evacuation for $\leq^{op}_{KL}$ and
$\leq_{weak}$, follow from the fact that the involutive maps
$$
 w \mapsto w_0w ~\text{ and }~w \mapsto ww_0
 $$
are antiautomorphisms of both $(\sym_n,\leq^{op}_{KL})$
\cite{Garsia-McLarnan} and  $(\sym_n,\leq_{weak})$.  Hence $w
\mapsto w_0ww_0$ is an automorphism  of  both. On the other hand
$P(ww_0)$ is just the  transpose tableau of $P(w)$ \cite{Schensted}
and $P(w_0ww_0)$ is nothing but the evacuation of $P(w)$
\cite{Schutzenberger1}.

Indeed  $w_0w$ and $ww_0$ correspond reversing the value and the
order of numbers in $w$ respectively. Therefore by  Greene's theorem
they reverse the dominance order on the $RSK$ insertion shapes which
then gives the desired property for $(SYT_n,\leq_{chain})$.

The assertion that Sch\"utzenberger's evacuation map gives a poset
automorphism of $(SYT_n,\leq_{geom})$ follows from Melnikov's work
\cite[Page 17--18]{Melnikov2}.
\end{proof}

\vskip .1in
\begin{question} (See discussion by Van Leeuwen \cite[\S8]{Leeuwen1}). Is  the
map
 which sends a tableau to its transpose an anti-automorphism of  the
 geometric order?
\end{question}

\vskip .1in
 By part (ii)  of the Proposition
\ref{order-preserving-maps}, if $S\leq T$ under the weak, $KL$,
geometric or chain orders then $\sh(S)\leq^{op}_{dom} \sh(T) $.
Actually we have a stronger condition for the first three orders
which is given in Proposition \ref{change-of-shapes} below. On the
other hand Example \ref{chain-order-shape-agree} shows that this
property is not satisfied  by chain order.

\begin{proposition}
\label{change-of-shapes}
Let $\leq$ be any of $\leq_{weak}$,
$\leq^{op}_{KL}$ or $\leq_{geom}$ on $SYT_n$. Then
 $$
 S \lneqq T ~~\Longrightarrow~~ \sh(S) \lneqq^{op}_{dom}\sh(T)
 $$
$e. g.$, under these orders  the shape of the tableaux is not fixed.
\end{proposition}

\begin{proof} For $\leq^{op}_{KL}$, this property can be induced from
the work of  Lusztig \cite{Lusztig}  which result in the conclusion
that, for $\sym_n$ right cells given by the tableaux of the same
shape form an antichain in the $KL$ order.

 For $\leq_{geom}$,  Gerstenhaber's result
mentioned in the proof of Proposition
\ref{order-preserving-maps}(ii) gives the required property; if
$\sh(S)=\sh(T)=\lambda$, the orbital varieties  $\mathcal{V}_T$ and
$\mathcal{V}_S$ lie in the same nilpotent orbit
$\mathcal{O}_\lambda$. As being the irreducible components of
$\mathcal{O}_\lambda \cap \mathfrak{n}$ they satisfy neither
$\mathcal{V}_T \subseteq \overline{\mathcal{V}}_S$ nor
$\mathcal{V}_S \subseteq \overline{\mathcal{V}}_T$. Therefore $T$
and $S$ are not comparable under $\leq_{geom}$ and this proves the
hypothesis.

Now $\leq_{weak}$  satisfy the hypothesis since it is weaker then
$KL$ and geometric orders.
\end{proof}

\begin{example}
\label{chain-order-shape-agree} The following tableaux have
$T\gneqq_{chain}S$   although they have the same shape.
$$
T=
\begin{array}{ccc} 1&3&6\\2&4\\5\end{array}
\hskip .1in \text{and} \hskip .1in S=\begin{array}{ccc}
1&3&4\\2&6\\5
\end{array}
$$
\end{example}

\vskip .15in
\subsection{Embedding}
\noindent

It is known  that  the (right) weak order on $\sym_n$ is weaker than
the (right) $KL$ preorder  on $\sym_n$ \cite[page
171]{Kazhdan-Lusztig}. As it is described, for instance in
\cite[page 9]{Melnikov1}, the weak order is also weaker than
geometric order  on $\sym_n$.   Therefore by the virtue of its
definition $(SYT_n,\leq_{weak})$ embeds in  $(SYT_n,\leq_{KL})$ and
$(SYT_n,\leq_{geom})$.

On the other hand by Corollary~\ref{segment-restriction-lemma} and
by Proposition~\ref{order-preserving-maps}$\mathrm{(ii)}$  the weak,
$KL$ and geometric orders on $SYT_n$ are weaker then chain order.

The following important result, which reveals that $KL$  order
embeds in geometric order on $SYT_n$, can be deduced from the work
of Melnikov \cite[Corollary 1.2]{Melnikov3}, Borho and Brylinski
\cite[6.3]{Borho} and Vogan \cite{Vogan}.

\begin{theorem} \label{KL-weaker-than-geometric}
On $\sym_n$, $KL$ order is weaker than geometric order. Therefore
for all $S, T \in SYT_n$,
$$
 S \leq^{op}_{KL} T ~~\Longrightarrow~~ S\leq_{geom} T.
 $$
\end{theorem}

It happens that all these  four orders coincide for $n \leq 5$, but
they start to differ for $n=6$.  Proposition \ref{change-of-shapes}
and the Example \ref{chain-order-shape-agree} provided above show
that $(SYT_n, \leq_{chain})$  differs from all the other orders for
$n=6$.  The following examples reveals the same fact for $(SYT_n,
\leq_{weak})$. (cf. Melnikov \cite[Example 4.1.6]{Melnikov1}).

\begin{example}
\label{KL-stronger-example}
 Let $
  S=\begin{array}{ccc} 1&2&3\\4&5&6\end{array},\hskip .1in
  T_1=\begin{array}{ccc} 1 &2& 5 \\ 3 &6 & \\ 4 & &
\end{array}
~\text{ and }~ T_2=\begin{array}{ccc} 1 &3& 6 \\ 2 &4&
\\ 5 &
\end{array}.
$

Computer calculations show that $S \leq^{op}_{KL} T_1, T_2$, but $S
\not\leq_{weak} T_1, T_2$. By using the anti-automorphism of
$\leq^{op}_{KL}$ and $ \leq_{weak}$ that transposes a standard Young
tableau (see Proposition~\ref{order-preserving-maps}) one obtains
two more examples of pairs of tableaux which are comparable in
$\leq^{op}_{KL}$, but not in $\leq_{weak}$.  These are the {\it
only} such examples in $SYT_6$.
\end{example}

To summarize we have the following diagram:

$$ \label{embedding-sequence} (SYT_n,\leq_{weak}) ~\subsetneqq~
(SYT_n,\leq^{op}_{KL}) ~\subseteq~ (SYT_n,\leq_{geom}) ~\subsetneqq~
(SYT_n,\leq_{chain}).
$$

\vskip .2in
\begin{question}
 Do $(SYT_n, \leq_{KL}^{op})$ and $(SYT_n, \leq_{geom})$ coincide?
\end{question}

\vskip .3in

\subsection{Extension  from  segments}
\label{Embedding-from-initial-segments-section} \noindent

In this section  we discuss  two  order preserving maps  which embed
$SYT_n$ into $SYT_{n+1}$ under any of the four orders.

 Denoted by $\Omega_1$ and
$\Omega_2$, these maps are given by the following rule: For each $T
\in SYT_n$, $\Omega_1:SYT_n \mapsto SYT_{n+1}$ concatenates $n+1$ to
the first row of $T$  from  the right whereas $\Omega_2:SYT_n
\mapsto SYT_{n+1}$ concatenates $n+1$ to the first column  of $T$
from  the bottom i.e.,

$\put(0,70){\makebox(1,1){}}
 \put(40,60){\line(1,0){50}}
 \put(70,40){\line(1,0){20}}
\put(40,20){\line(1,0){30}} \put(40,20){\line(0,1){40}}
 \put(70,20){\line(0,1){20}}
\put(90,40){\line(0,1){20}}
\put(110,20){\makebox(12,35){$\underset{\Omega_1}{\longrightarrow}$}}
\put(140,60){\line(1,0){50}}
 \put(170,40){\line(1,0){20}}
\put(140,20){\line(1,0){30}} \put(140,20){\line(0,1){40}}
 \put(170,20){\line(0,1){20}}
\put(190,40){\line(0,1){20}}
\put(190.5,49.5){\framebox(20,10){{\scriptsize $n+1$}}}
\put(260,60){\line(1,0){50}}
 \put(290,40){\line(1,0){20}}
\put(260,20){\line(1,0){30}} \put(260,20){\line(0,1){40}}
 \put(290,20){\line(0,1){20}}
\put(310,40){\line(0,1){20}}
\put(330,20){\makebox(12,35){$\underset{\Omega_2}{\longrightarrow}$}}
\put(360,60){\line(1,0){50}}
 \put(390,40){\line(1,0){20}}
\put(360,20){\line(1,0){30}} \put(360,20){\line(0,1){40}}
 \put(390,20){\line(0,1){20}}
\put(410,40){\line(0,1){20}}
\put(360,10){\framebox(20,10){{\scriptsize $n+1$}}} $

\begin{definition} Any partial order $\leq$  on $SYT_n$ is said to have  the
property of {\it
extension from  segments} if the maps $\Omega_1, \Omega_2: SYT_n
\mapsto SYT_{n+1}$
 are order preserving.
\end{definition}

In what follows we will prove that all of the four orders have the
extension from segments property.

\begin{lemma}
\label{embedding-from-initial-segments-theorem}The  maps  $\Omega_1$
and $\Omega_2$ are order preserving under  the weak, $KL$, the
geometric and chain orders.
\end{lemma}

\begin{proof}
For any $T \in SYT_n$ and $\tau \in C_T$, let $\tau(n+1)$ and
$(n+1)\tau$ be the words obtained by concatenating $n+1$ to $\tau$
from the right and respectively from the left. The $RSK$ insertion
algorithm yields that
$$P(\tau(n+1))=\Omega_1(T) ~\text{ and }~ P((n+1)\tau)=\Omega_2(T).$$
Conventionally, we use the following notation:
$$\Omega_1(\tau):=\tau(n+1) ~\text{ and }~ \Omega_2(\tau):=(n+1)\tau.
$$

\vskip .1in
 \noindent {\it Chain order}: Let $S \leq_{chain} T$  in
$SYT_n$, i.e., for any $1\leq i<j\leq n$ one has
$$\sh(\st(S_{[i,j]})) \leq_{dom}^{op} \sh(\st(T_{[i,j]})).$$
Now concatenating  $n+1$ to the first row of $S$ and $T$ from the
right (after applying jeu de taquin slides)  obviously does not
affect $\sh(\st(S_{[i,j]}))$ and $\sh(\st(T_{[i,j]}))$ if $j< n+1$,
and both have $n+1$ added to first row  if $j=n+1$. Therefore
$$\Omega_1(S) \leq_{chain} \Omega_1(T).$$

On the other hand by
Proposition~\ref{order-preserving-maps}($\mathrm{iv}$) one has:
$$S~\leq_{chain}~T ~\Longrightarrow~S^t~\geq_{chain}~T^t ~\Longrightarrow~
\Omega_1(S^t)^t ~\leq_{chain}~ \Omega_1(T^t)^t$$
 and since $\Omega_1(S^t)^t=\Omega_2(S)$ for any tableau $S$, now $\Omega_2$ is
also order preserving.

\vskip .1in
 \noindent {\it Weak order}: For this it is enough to
consider the covering relations of $(SYT_n, \leq_{weak})$. If $S$ is
covered by $T$ then there exist two permutations $\sigma \in C_S$
and $\tau \in C_T$ such that $ \sigma \leq_{weak} \tau$.
Equivalently $\Inv_{L}(\sigma) \subset \Inv_{L}(\tau)$. On the other
hand the last assertion implies
$$\Inv_{L}(\Omega_1(\sigma)) \subset \Inv_{L}(\Omega_1(\tau)) ~\text{
and }~\Inv_{L}(\Omega_2(\sigma)) \subset \Inv_{L}(\Omega_2(\tau)).$$
Therefore in either case the weak order relation is preserved and we
have
$$\Omega_1(S) \leq_{weak}\Omega_1(T) ~\text{ and }~ \Omega_2(S) \leq_{weak}
\Omega_(T).$$

\vskip .1in \noindent
{\it $KL$ order}:  This fact for  $KL$ order
can be deduced easily by considering $\sym_n$ as a parabolic
subgroup of $\sym_{n+1}$: any two permutations $v, w \in \sym_n$
satisfying  $v \leq_{KL}^{op} w$ in the parabolic subgroup  $\sym_n$
still have the same relation in $\sym_{n+1}$.

If $S\leq_{KL}^{op} T$ then  there exist $\sigma \in C_S$ and $\tau
\in C_T$ satisfying $\sigma \leq_{KL}^{op} \tau$ in $\sym_n$.  Then
concatenating $n+1$ to the right side of both words still yields
$\Omega_1(\sigma)  \leq_{KL}^{op} \Omega_1(\tau) $ in $\sym_{n+1}$.
Hence $\Omega_{1}(S) \leq_{KL}^{op} \Omega_1(T)$ and by
Proposition~\ref{order-preserving-maps}($\mathrm{iv}$)
$\Omega_{2}(S) \leq_{KL}^{op} \Omega_2(T)$.

\vskip .1in \noindent
{\it Geometric order}: This fact follows from
the result of Melnikov \cite[Proposition 6.6]{Melnikov4}.

\end{proof}

\vskip  .2in

\subsection{Extension by $RSK$ insertions}
\label{extension-by-RSK-insertions-section} \noindent

In  \cite{Melnikov1},  Melnikov indicates another extension property
of the weak and geometric order on $SYT_n$ which also generalize the
property of extension  from segments. \vskip .1in Let $\leq$~ be any
order on $SYT_n$, ~$i\leq n$ and   $\acute{S}$ and $\acute{T}$ are
some tableaux on $[n]-\{i\}$. Suppose  $S$ and $T$ are the tableaux
in $SYT_{n-1}$ obtained by standardizing $\acute{S}$ and $\acute{T}$
respectively. Define  an order on $\acute{S}$ and $\acute{T}$ in the
following way
$$
\acute{S} \leq \acute{T} \hskip .1in \text{ if } \hskip .1in S \leq
T.
$$
Then   ~$\leq$~ is said to have   the property of {\it extension by
$RSK$ insertions} if the $RSK$ insertion of the  element $i$ into
both tableaux $\acute{S}$ and $\acute{T}$ from above (or from the
left) still preserves the order, in other words, denoting the
resulting tableaux by $\acute{S}^{\downarrow i}$ and
$\acute{T}^{\downarrow i}$, if one  has
$$
\acute{S} \leq \acute{T} ~\Rightarrow~\acute{S}^{\downarrow i} \leq
\acute{T}^{\downarrow i}.
$$

The property of extension by $RSK$ insertions  for the weak order
and geometric order was first proven by Melnikov in \cite{Melnikov1}
and \cite{Melnikov4} respectively. The same fact for $KL$ order can
be deduced from the work of Barbash and Vogan \cite[2.34,
3.7]{Barbash-Vogan} by using the theory that relates Kazhdan-Lusztig
(left) cells  to  primitive ideals. Below, independently from this
theory, we provide a proof that shows $KL$ order has the property of
extension by $RSK$ insertions.  On the other hand the following
example shows that chain order does not have this property.
$$
\acute{T}=\begin{array}{ccc} 1&3&7\\2&5\\6 \end{array} \geq_{chain}
\begin{array}{ccc} 1&3&5\\2&7\\6 \end{array}=\acute{S} \hskip .2in
\text{but} \hskip .2in \acute{T}^{\downarrow 4}=\begin{array}{ccc}
1&3&4\\2&5&7\\6 \end{array} \not \geq_{chain}~
\begin{array}{ccc} 1&3&4\\2&5\\6&7 \end{array}
=\acute{S}^{\downarrow 4}.
$$

\begin{lemma}\label{KL-RSK-insertions-property Lemma}
$KL$ order on $SYT_n$ has the extension by $RSK$ insertions
property.
\end{lemma}
\begin{proof}
Let $\acute{S}$ and $\acute{T}$ be two tableaux on $[n]-\{i\}$ such
that  $\acute{S} \leq_{KL}^{op} \acute{T}$. In other words for $S$
and $T$ which are obtained by standardizing  $\acute{S}$ and
$\acute{T}$ respectively, we have $S \leq_{KL}^{op} T $. We may
assume that $S$  is covered by $T$. Then there exist $\sigma$ and
$\tau$ in the Knuth classes of $S$ and $T$ respectively such that $
\sigma ~\lessdot_{KL}^{op}~ \tau ~\text{ in }~\sym_{n-1}$. Since
$\sym_{n-1}$ is a parabolic subgroup of $\sym_n$,  as
Lemma~\ref{embedding-from-initial-segments-theorem} for the $KL$
order shows, concatenating $n$ to the right side of both
permutations yields $\sigma n ~\lessdot_{KL}^{op} ~\tau n$ in
$\sym_n$. Therefore we have
$$
D_R(\tau n)- D_R(\sigma n) \not = \varnothing ~\text{ and }~
\begin{cases}
\text{ either }& \sigma n \leq \tau n
~\text{ and}~\overline{\mu}(\sigma n,\tau n)\not=0 \\
\text{ or }&  \tau n\leq \sigma n \text{ and } ~\overline{\mu}(\tau
n,\sigma n)\not=0
\end{cases}
$$
where $\leq$ denotes Bruhat order. Without lost of generality we
assume ~$\sigma n \leq \tau n$ ~and ~$\overline{\mu}(\sigma n,\tau
n)\not=0$.

\vskip .1in
 Consider the  permutations~ ~$s_{i}s_{i+1}\ldots s_{n-1}(\sigma n) ~\text{ and
}~s_{i}s_{i+1}\ldots s_{n-1}(\tau n)$~ ~which are obtained by
multiplying $\sigma n$ and $\tau n$ from the left by the
transpositions ~$s_{n-1}, s_{n-2}, \ldots, s_{i+1}, s_{i}$ ~in this
order.  It is easy to check that the $RSK$ insertion tableaux of
~$s_{i}s_{i+1}\ldots s_{n-1}(\sigma n)$~ and ~$s_{i}s_{i+1}\ldots
s_{n-1}(\tau n)$~ are nothing but $\acute{S}^{\downarrow i}$ and
respectively $\acute{T}^{\downarrow i}$. Then $\acute{S}^{\downarrow
i} \leq_{KL}^{op} \acute{T}^{\downarrow i}$ follows, once it is
shown that
\begin{equation}\label{extension-by-RSK-insertions-KL-eq}
 s_{i}s_{i+1}\ldots s_{n-1}(\sigma n) ~\leq_{KL}^{op}
~s_{i}s_{i+1}\ldots s_{n-1}(\tau n).
\end{equation}

 Let
 $$u_n=\sigma n \text{ and } w_n=\tau n
 $$
 and for each  $k$ such that $i\leq k\leq n-1$, let
$$
u_k=s_{k}\ldots s_{n-1}(\sigma n)~~\text{ and  }~~ w_k=s_{k}\ldots
s_{n-1}(\tau n).
$$

 Obviously  for each $i \leq k \leq n$, analysis on the (left) inversion sets
yields
\begin{equation}\label{extension-by-RSK-insertions-KL-eq2}
l(u_k,w_k)=l(\sigma n,\tau n)
\end{equation}
and one can check that
\begin{equation}\label{extension-by-RSK-insertions-KL-eq3}
u_k \leq w_k
\end{equation}
by using a basic characterization of Bruhat order. That is: $u \leq
w$ in $\sym_n$ if and only if for each $j\leq n$, the sets of the
form $ \{u_1,\ldots,u_j\}  ~\text{ and }~\{w_1,\ldots,w_j\}$ can be
compared in the manner that after ordering their elements from the
smallest to the biggest, the $i$-th element of the first set is
always smaller than or equal to the $i$-th element of the second set
for each $i\leq j$.

On the other hand
 multiplying $\sigma n$~ and $\tau n$~ by ~$s_{k} \ldots
s_{n-2} s_{n-1}$ from the left does not change the right descents of
these permutations on the first $n-1$ positions. In other words,
when  restricted to the first $n-1$ positions $\sigma n$ and $u_k$
(similarly $\tau n$ and $w_k$) share the same right descents.
Therefore
\begin{equation}\label{extension-by-RSK-insertions-KL-eq4}
 D_R(\tau n)- D_R(\sigma n)
\not = \varnothing \Longrightarrow D_R(w_k)- D_R(u_k) \not =
\varnothing.
\end{equation}

 Now we will show  that
$$
P_{u_k,w_k}(q)=P_{\sigma n,\tau
n}(q).
$$
 Obviously $P_{u_n,w_n}(q)=P_{\sigma n,\tau n}(q)$ and
therefore  it is enough to prove that
$P_{u_k,w_k}(q)=P_{u_{k+1},w_{k+1}}(q)$, since then the required
equality follows by induction. \vskip .05in Observe that $u_k=s_k
u_{k+1}$, $w_k=s_k w_{k+1}$ i.e.,  both of them are permutations in
$\sym_n$ ending with the number $k$. So $s_k$ lies both in
$D_L(u_k)$ and $D_L(w_k)$ and  by
\eqref{recursive-formula-KL-polynomial}
$$
P_{u_k,w_k}(q)=P_{u_{k+1},w_{k+1}}(q)+q P_{u_k,w_{k+1}}(q)-
\sum_{\{v:s_k\in D_L(v)\}} q^{l(v,w)/2}~
\overline{\mu}(v,w_{k+1})~P_{u_k,v}(q).
$$
Since $u_k$ ends with $k$ and  $w_{k+1}$ ends with $k+1$, from the
characterization of the Bruhat order it follows that ~$u_k \not\leq
w_{k+1}$ and furthermore  there exist no permutation $v$ satisfying
$u_k \leq v \leq w_{k+1}$. Then by
\eqref{basic-properties-KL-polynomial}, all the summation terms on
the right hand side, except $P_{u_{k+1},w_{k+1}}(q)$, are equal to
$0$. Henceforth
$$
P_{u_k,w_k}(q)=P_{u_{k+1},w_{k+1}}(q)
$$
and $P_{u_k,w_k}(q)=P_{\sigma n, \tau n}(q)$ follows by induction.
This result together with \eqref{KL-Mu-function} and
\eqref{extension-by-RSK-insertions-KL-eq2} imply that
\begin{equation}\label{extension-by-RSK-insertions-KL-eq5}
\overline{\mu}(u_k,w_{k})=\overline{\mu}(\sigma n, \tau n)\not = 0.
\end{equation}
 Therefore  by
\eqref{KL-preorder}, \eqref{extension-by-RSK-insertions-KL-eq3},
\eqref{extension-by-RSK-insertions-KL-eq4} and
\eqref{extension-by-RSK-insertions-KL-eq5} we have  $u_k
\leq_{KL}^{op} w_k$ for eack $i \leq k \leq n-1$ and so
\eqref{extension-by-RSK-insertions-KL-eq} is true. Hence
$\acute{S}^{\downarrow i} \leq_{KL}^{op} \acute{T}^{\downarrow i}.$
\end{proof}

\vskip .3in
\section{Proof of Theorem \ref{P-R-theorem}}
\label{P-R-theorem-section}
\vskip .1in

Malvenuto and Reutenauer, in \cite{Malvenuto-Reutenauer} construct
two graded Hopf algebra structures on the  $\mathbb{Z}$ module of
all permutations $\mathbb{Z}S = \oplus_{n \geq 0} \mathbb{Z}S_{n}$
which are dual to each other, and shown to be free as associative
algebras by Poirier and Reutenauer in \cite{Poirier-Reutenuer}. The
product structure of the one that concerns us here is given by
$$
u \ast w:= \shf(u,\overline{w})
$$
where $\overline{w}$ is obtained by increasing the indices of $w$ by
the length of $u$ and $\shf$ denotes the shuffle product.
\vskip .1in
Poirier and Reutenauer also show that $\mathbb{Z}$ module of all
plactic classes $\{PC_{T}\}_{T \in SYT}$, where
$PC_{T}=\sum_{\footnotesize P(u)=T}u$ becomes a Hopf subalgebra of
permutations whose product (also shown  in \cite{Lascoux} and
\cite{Sottile}) is given by the formula
\begin{equation} \label{P-R-multiplication}
PC_{T} \ast PC_{T'} =
        \sum_{\substack{P(u)=T \\  P(w)=T'}}
    \shf(u, \overline{w }) \end{equation}

Then the bijection sending each plactic class to its defining
tableau gives us a Hopf algebra structure on the $\mathbb{Z}$ module
of all standard Young tableaux, $\mathbb{Z}SYT = \oplus_{n \geq 0}
\mathbb{Z}SYT_{n}$.

\vskip .1in For example,
\begin{equation}\label{P-R-example}
\begin{array}{ll}
\begin{array}{l}
PC_{\footnotesize{\begin{array}{ll} 1&2 \\ 3&
\end{array}}}~\ast~
    PC_{\footnotesize{\begin{array}{l} 1\\ 2  \end{array}}}\\
    \\ \\ \\ \\ \\ \\
    \end{array}
&
\begin{array}{ll}
=& \shf(312,54) + \shf(132,54)\\
=& 31254 + 31524 +35124 + 53124+31542+ 35142+ 35412+ \\
&53142 +  53412+54312 +13254 + 13524 +15324 + 51324+\\
&13542+ 15342+15432+ 51342+51432+54132\\ \\=&
PC_{\footnotesize{\begin{array}{lll} 1&2&4\\3&5
\end{array}}}+
 PC_{\footnotesize{\begin{array}{lll} 1&2&4\\3\\ 5 \end{array}}}+
 PC_{\footnotesize{\begin{array}{ll} 1&2\\3&4\\ 5 \end{array}}}+
 PC_{\footnotesize{\begin{array}{ll} 1&2\\3\\4\\5 \end{array}}}.
\end{array}
\end{array}
\end{equation}

Another approach to calculate the product of two tableaux is given
in \cite{Poirier-Reutenuer} where Poirier and Reutenauer explain
this product using jeu de taquin slides. Our goal is to show that it
can also be described by a formula using partial orders, analogous
to a result of Loday and Ronco \cite[Thm. 4.1]{Loday-Ronco}. To
state their result, given $\sigma \in \sym_k$ and $\tau \in
\sym_\ell$, with $n:=k + \ell$, let $\overline{\tau}$ be obtained
from $\tau$ by adding $k$ to each letter.  Then let $\sigma/\tau$
and $\sigma\backslash\tau$ denote the concatenations of $\sigma,
\overline{\tau}$ and of $\overline{\tau}, \sigma$, respectively.

\begin{theorem}
\label{Loday-Ronco-Malvenuto-Reutenauer} For $\tau \in \sym_k$ and
$\sigma \in \sym_\ell$, with $n:=k + \ell$, one has in the
Malvenuto-Reutenauer Hopf algebra
$$
\sigma \ast \tau =
 \sum_{\substack{\rho \in \sym_n:\\
      \sigma/\tau \leq \rho \leq \sigma\backslash \tau}} \rho.
$$
Equivalently, the shuffles $\shf(\sigma,\overline{\tau})$ are the
interval $[\sigma/\tau , \sigma\backslash \tau]_{\leq_{weak}}.$
\end{theorem}

\vskip .1in
The following facts  are crucial  for transporting the Loday and
Ronco result to $SYT_n$.

Let $\sigma \in \sym_k, \tau \in \sym_\ell$. When $P(\sigma)=S$ and
$P(\tau)=T$, let $\overline{T}$ denote the result of adding $k$ to
every entry of $T$.  It is easily seen that
\begin{equation}\label{P-R-sileds}
 P(\sigma/\tau) = S/T \mbox{ } \mathrm{ and } \mbox{ } P(\sigma\backslash\tau) =
S\backslash T
\end{equation}
where $S/T$ (respectively, $S\backslash T$) is the tableaux whose
columns (resp. rows) are obtained by concatenating the columns
(resp. rows) of $S$ and $\overline{T}$. Note also that
Lemma~\ref{j-d-t-initial-final} shows for $I=[k]$
\begin{equation}\label{P-R-restriction}
\begin{aligned}
(S/T)_I&=S & \st((S/T)_{I^c})&=T\\
(S\backslash T)_I&=S & \st((S\backslash T)_{I^c})&=T.
\end{aligned}
\end{equation}

\vskip .2in

 The following  theorem is a consequence of Lemma
 \ref{j-d-t-initial-final},
 Corollary \ref{segment-restriction-lemma} and Theorem
\ref{Loday-Ronco-Malvenuto-Reutenauer}.

\begin{theorem}
\label{P-R-theorem-0}  Let ~$\leq$~ be a partial order on $SYT_n$,
for all $n>0$, that
\begin{enumerate}
\item[(a)] is stronger than  $\leq_{weak}$ and
\item[(b)] restricts to segments.
\end{enumerate}
Then in the Poirier-Reutenauer Hopf algebra,
 $$ S \ast T = \sum_{\substack{R \in SYT_n:\\
                 S/T \leq R \leq S\backslash T}} R. $$
\end{theorem}

\begin{proof}
Let $\leq$ be a partial order on $STY_n$ satisfying hypothesis.
 From (\ref{P-R-multiplication}), (\ref{P-R-sileds}) and
Theorem \ref{Loday-Ronco-Malvenuto-Reutenauer} it follows that any
tableau $R$ appearing in the product $ S \ast T$ satisfies:  $S/T
\leq_{weak} R \leq_{weak} S\backslash T$. Therefore  we have  $S/T
\leq R \leq S\backslash T$ and this proves one direction.

Let $R$ be any tableau such that $S/T \leq R \leq S\backslash T$.
Also let  $I=[k]$ where $k$ is the size of the tableau $S$ .  By
hypothesis
$$S=(S/T)_I \leq  R_I  \leq (S\backslash T)_I=S $$
$$T=\st((S/T)_{I^c})\leq \st(R_{I^c}) \leq \st((S\backslash T)_{I^c})=T$$
i.e.,   $R_I=S$ and  $\st(R_{I^c})=T$ and this shows that $R$ can be
found by shuffling  $S$ and $T$ in a certain way. Therefore $R$ lies
in the product $ S \ast T$.
\end{proof}

\begin{proof} [Proof of Theorem~\ref{P-R-theorem}]
All four orders
$\leq_{chain}$, $\leq_{weak}$, $\leq_{KL}^{op}$ and $\leq_{geom}$ on
$SYT_n$ satisfy the hypotheses of Theorem~\ref{P-R-theorem-0} by
Corollary~\ref{segment-restriction-lemma}. Therefore the result
follows.
\end{proof}

\begin{example}
Let $ T= \scriptstyle{\begin{array}{cc} 1&2 \\ 3 & \end{array}}$ and
  $S =\scriptstyle{\begin{array}{c} 1 \\ 2 \end{array}}$.
Then  $T/S=\scriptstyle{\begin{array}{ccc} 1&2&4 \\ 3 & 5
\end{array}}$, ~ $T\backslash S=\scriptstyle{\begin{array}{cc} 1&2 \\ 3& \\ 4&
\\ 5
\end{array}}$ and
  \eqref{P-R-example} gives
$$T \ast S= \begin{array}{cc} 1&2 \\ 3 \end{array} \ast
\begin{array}{c} 1 \\ 2  \end{array}=
\begin{array}{ccc} 1 & 2 &4  \\3& 5  \end{array} +
\begin{array}{ccc} 1 & 2 &4  \\3 \\ 5  \end{array} +
\begin{array}{cc} 1 & 2   \\3& 4 \\5  \end{array} +
\begin{array}{cc} 1& 2 \\ 3 \\  4 \\ 5 \end{array}.$$
On the other hand, when considered with any of the four orders, the
Hasse diagram of  $SYT_{5}$ in Figure~\ref{figure1} shows that the
product above is equal to the sum of all tableaux in the interval
$[T/S, T\backslash S]$.
\end{example}

\vskip .3in
\section{Proof of Theorem~\ref{homotopy-theorem}}
\label{homotopy-theorem-section}
 \vskip .1in

 In this section, we prove
Theorem~\ref{homotopy-theorem}.  We will view the commutative
diagram
\begin{equation}
\label{SYT-triangle}
\begin{array}{cll}
        S_{n} & \longrightarrow & SYT_{n} \\

         &  \searrow   & \downarrow  \\
                           &       &    2^{[n-1]} \\
\end{array}
\end{equation}
as an instance of the following set-up, involving closure relations,
equivalence relations, order-preserving maps, and the topology of
posets.  For background on poset topology, see \cite{Bjorner}.

Let $P$ be a partially  ordered set $(P, \leq_P)$  and $p \mapsto
\bar{p}$ a closure relation on $P$, that is,
$$\begin{aligned}
\bar{\bar{p}}&=\bar{p}, &
p &\leq_P \bar{p}& \text{and }
p &\leq_P q \text{ implies }\bar{p} \leq_P \bar{q}.
\end{aligned}$$
It is well-known  \cite[Corollary 10.12]{Bjorner} that in this
instance, the order-preserving closure map $P \rightarrow
\overline{P}$ has the property that its associated simplicial map of
order complexes $\Delta(P) \rightarrow \Delta(\overline{P})$ is a
strong deformation retraction.

Now assume $\sim$ is an equivalence relation on $P$ such that, as
maps of sets, the closure map $P \rightarrow \overline{P}$ factors
through the quotient map $P \rightarrow P/_{\sim}$.  Equivalently,
the vertical map below is well-defined, and makes the diagram
commute:

\begin{equation}
\label{commuting-set-maps}
\begin{array}{cll}
        P & \longrightarrow & P/_{\sim} \\
        &  \searrow   & \downarrow  \\
        &       &    \overline{P} \\ \end{array}
\end{equation}

\begin{proposition}
\label{closure-factorization} In the above situation, partially
order $\overline{P}$ by the restriction of $\leq_P$, and assume that
$P/_{\sim}$ has been given a partial order $\leq$ in such a way that
the horizontal and vertical maps in the \eqref{commuting-set-maps}
are also order-preserving. Then the commutative diagram of
associated simplicial maps of order complexes are all homotopy
equivalences.
\end{proposition}

\begin{proof}
Obviously one can define a closure relation on
$P/_{\sim}$ such that $\overline{P/_{\sim}}=\overline{P}$, and the
result follows.
\end{proof}

\begin{lemma}
\label{descent-fibers}
Given any subset $D \subset [n-1]$, there exists a maximum element
$\tau(D)$ in $(\sym_n, \leq_{weak})$ for the descent class
$$\Des_L^{-1}(D):=\{ \sigma \in \sym_n: \Des_L(\sigma) = D\}.$$
Consequently, the map $\sym_n \rightarrow \sym_n$ defined by $\sigma
\mapsto \tau(\Des_L(\sigma))$ is a closure relation which also
restricts to the proper parts and  its  image is isomorphic to
$(2^{[n-1]},\subseteq)$.
\end{lemma}

\begin{proof}
 It is known that  \cite[page 98-100]{Bjorner2}
$$\Des_L^{-1}(D):=\{ \sigma \in \sym_n: \Des_L(\sigma) = D\}$$
 is actually an interval of the weak Bruhat order on $S_{n}$.  Therefore
 the map $\sigma \mapsto \tau(\Des_L(\sigma))$ is a closure relation
 and since  $\Des_L^{-1}(\varnothing)$ and
$\Des_L^{-1}([n-1])$ consist of respectively $\hat{0}$  and
$\hat{1}$ in ($\sym_n,\leq_{weak}$), it restricts to the proper
parts. Now it is easy to see that its image
 is isomorphic to $(2^{[n-1]},\subseteq)$.
\end{proof}

\begin{corollary}
\label{SYT-homotopy-corollary} Order $\sym_n$ by $\leq_{weak}$ and
$2^{[n-1]}$ by $\subseteq$. Let $\leq$ be any order on $SYT_n$ such
that the commuting diagram \eqref{SYT-triangle} has all the maps
order-preserving.  Then these restrict to a commuting diagram of
order-preserving maps on the proper parts, each of which induces a
homotopy equivalence of order complexes.  Consequently,
$\mu(\hat{0},\hat{1}) = (-1)^{n-3}$.
\end{corollary}

\begin{proof}
The fact that the maps restrict to the proper parts follows because
we know the maps explicitly as maps of sets, and the images of
$\hat{0}, \hat{1}$ in $(\sym_n, \leq_{weak})$  must be exactly the
$\hat{0}, \hat{1}$ in $(SYT_n, \leq)$ (namely the single-row and
single-column tableaux) because the horizontal map is
order-preserving.

The fact that they induce homotopy equivalences follows from
Proposition~\ref{closure-factorization} applied to the three proper
parts, using the closure relation in Lemma~\ref{descent-fibers} and
letting $\sim$ be Knuth equivalence $\Knuth$.  One must observe that
$\Des_L(\sigma)$ depends only on the Knuth class of $\sigma$.

The fact that $\mu(\hat{0},\hat{1}) = (-1)^{n-3}$ for the Boolean
algebra $(2^{[n-1]}, \subseteq)$ is well-known
\cite[Prop.3.8.4]{Stanley}.
\end{proof}

\vskip .1in
\begin{proof}[Proof of Theorem~\ref{homotopy-theorem}] By
Proposition~\ref{order-preserving-maps} all four orders   on $SYT_n$
 satisfy the hypotheses of
Corollary~\ref{SYT-homotopy-corollary}.
\end{proof}

\vskip .1in
\begin{example}
In Figure~\ref{mob}, the first interval in $(SYT_8,\leq_{chain})$
has M\"obius function value $2$, whereas  M\"obius function value of
the second interval which is found in $(SYT_8,\leq_{weak})$,
$(SYT_8,\leq_{KL}^{op})$ and  $(SYT_8,\leq_{geom})$ is $-2$.
Therefore the M\"obius function values of the proper intervals in
$\leq_{chain}$, $\leq_{weak}$, $\leq_{KL}^{op}$ and  $\leq_{geom}$
on $SYT_n$ need not all lie in $\{\pm 1, 0\}$ as they do in
$(\sym_n,\leq_{weak})$.
\end{example}

\vskip .3in
\section{Inner translation and skew orders}
\label{inner-translation-skew-orders-section}
\vskip .1in

In this section we describe the {\it inner translation property } of
$KL$ and geometric order on $SYT_n$ which enable us to generalize
these orders  to the skew standard Young tableaux $SYT_n^\mu$ of
size $n$ with some fixed inner boundary $\mu$.

To do this first we need to recall  the notion of {\it dual Knuth}
relations on $\sym_n$: permutations $\sigma, \tau \in \sym_n$ are
said to be {\it differ by a single dual Knuth relation} if for some
$i \in  [n-2] $, ~$i \in \Des_L(\sigma)$ and $i+1 \not \in
\Des_L(\sigma)$ whereas $i+1 \in \Des_L(\tau)$ and $i \not \in
\Des_L(\tau)$. In this case
$$\begin{aligned}
\mbox{ either } & \sigma=\ldots i+1 \ldots i \dots i+2 \ldots \mbox{
and }
\tau=\ldots i+2 \ldots i \dots i+1 \ldots \\
\mbox{ or } & \sigma=\ldots i+1 \ldots i+2 \dots i \ldots \mbox{ and
} \tau=\ldots i \ldots i+2 \dots i+1 \ldots
\end{aligned}$$

We say $\sigma, \tau $ are {\it Knuth equivalent} written as
$\sigma \nuth \tau $, if $\tau$ can be generated from $\sigma$ by a
sequence of single dual Knuth relations. Observe that $\sigma \nuth
\tau $ if and only if $\sigma^{-1} \Knuth \tau^{-1}$.

Since left descent sets are all equal for the  permutations in a
Knuth class $C_T$, $T \in SYT_n$, a single dual Knuth relation gives
the following action on tableaux: Let $r_T(i)$ be the row number of
$i$ in $T$ from the top.

 {\it Case 1}. If $i+1 \in \Des(T)$ and  $i
\not \in \Des(T)$ then
$$\begin{aligned}
\mbox{ either } && r_T(i+2)> r_T(i)\geq r_T(i+1)\\
\mbox{ or } && r_T(i) \geq r_T(i+2)> r_T(i+1).
\end{aligned}$$

The resulting tableau is found by interchanging $i+2$ and $i+1$ in
the first case and interchanging $i$ and $i+1$ in the second case.

{\it Case 2.} If $i \in \Des(T)$ and  $i+1 \not \in \Des(T)$ then
$$\begin{aligned}
\mbox{ either } && r_T(i+1)> r_T(i)\geq r_T(i+2)\\
\mbox{ or } &&  r_T(i+1) \geq r_T(i+2)> r_T(i) .
\end{aligned}$$

This time interchanging $i+2$ and $i+1$ in the first case and
interchanging $i$ and $i+1$ in the second case gives us  the
resulting tableau under the action of the  single dual Knuth
relation given with the triple $\{i,i+1,i+2\}$.

 We say $T\nuth T'$ if $T'$ can be obtained from
$T$ by applying a sequence of single dual Knuth relations as
described above. The following theorem, see \cite[Proposition
3.8.1]{Sagan} for example, is a nice characterization of this
relation.

\begin{theorem}
\label{dual-Knuth}
 Let $S, T \in SYT_n$. Then
$ S \nuth T ~~\text{ if and only if }~~ \sh(S)=\sh(T).$
\end{theorem}

\vskip .1in
Let $\{\alpha,\beta\}=\{i,i+1\}$ and
$SYT_n^{[\alpha,\beta]}$ be a subset of $SYT_n$ given by
$$SYT_n^{[\alpha,\beta]}:=\{ T \in SYT_n \mid
\alpha \in \Des(T), \beta \not \in \Des(T)\}.$$

As we described above we can apply  a single dual Knuth relation
determined with the triple $\{i, i+1, i+2\}$ on each  $ T\in
SYT_n^{[\alpha,\beta]}$ and this gives us the following  {\it inner
translation map}:
$$ \mathcal{V}_{[\alpha,\beta]}: SYT_n^{[\alpha,\beta]}\mapsto
SYT_n^{[\beta,\alpha]},$$ where $\mathcal{V}_{[\beta,\alpha]}\circ
\mathcal{V}_{[\alpha,\beta]}$ and $\mathcal{V}_{[\alpha,\beta]}\circ
\mathcal{V}_{[\beta,\alpha]}$ are just identity maps on their
domains.

\vskip .1in
\begin{definition} Any order $\leq$ on $SYT_n$ is said to have the {\it
inner translation property}  if the inner translation map
$$\mathcal{V}_{[\alpha,\beta]}: (SYT_n^{[\alpha,\beta]},\leq)\mapsto
(SYT_n^{[\beta,\alpha]},\leq)$$ is order preserving.
\end{definition}

\vskip .1in
 Now we give  the following corollary  which is crucial
in the sense that it provides the sufficient tool for  generalizing
any partial order on standard Young tableaux  to  the skew standard
tableaux.

\vskip .1in
 For  $1 \leq k<n$, let $R$ be a tableau in $SYT_k$ and
$$SYT_n^{R}:=\{ T  \in SYT_n \mid T_{[1,k]}=R\}.$$

\vskip .1in
\begin{corollary}
\label{equal-inner-shapes} Suppose $S, T \in SYT_n$  and $R,R' \in
SYT_k$ satisfy
$$S_{[1,k]}=T_{[1,k]}=R$$
$$\sh(R)=\sh(R').$$
Moreover  suppose  $S'$ and $T'$ are  the tableaux in $SYT_n$
obtained by replacing $R$ by $R'$ in $S$ and $T$ respectively.

Then for  $\leq$ having the inner translation property on $SYT_n$,
one has
  $$ S\leq T \mbox{ if and only if  } S'\leq T'.$$
In particular  $(SYT_n^{R},\leq)$ and $(SYT_n^{R'},\leq)$ are
isomorphic subposets of $(SYT_n,\leq)$.
\end{corollary}

\begin{proof} As a consequence of
 Theorem~\ref{dual-Knuth}, by applying to $S$ and $T$ the same  sequence of
dual Knuth relations on their  subtableau $R$, one can generate $S'$
and $T'$  respectively. On the other hand, since $\leq$ has inner
translation property at each step the order is preserved.
\end{proof}

\begin{theorem}
\label{Vogan-involution} $KL$ and  geometric order on $SYT_n$ have
the inner translation property. Therefore for any $R,R' \in SYT_k$
such that $\sh(R)=\sh(R')$ and $k<n$,
$$
(SYT_n^{R},\leq) ~\text{ and }~ (SYT_n^{R'},\leq)
$$
are isomorphic subposet of $SYT_n$ in  $KL$ and geometric orders.
\end{theorem}

\begin{proof}
This map is  first introduced by Vogan in ~\cite{Vogan} for  $KL$
order, where he also shows the desired property. For  geometric
order this result is  due to Melnikov ~\cite[Proposition
6.6]{Melnikov4}.
\end{proof}

The example given below shows that chain and the weak order do not
satisfy the inner translation property. (See also
Remark~\ref{inner-translation-remark}).

\begin{example}
$$\begin{array}{ccc} 1&3&\mathbf{6}\\2&\mathbf{4}\\\mathbf{5} \end{array}
\geq_{chain}
\begin{array}{ccc} 1&3&\mathbf{4}\\2&\mathbf{6}\\\mathbf{5}
\end{array}
\hskip .1in \text{ but } \hskip .1in
\begin{array}{ccc}
1&3&\mathbf{5}\\2&\mathbf{4}\\\mathbf{6}
\end{array} \not\geq_{chain}
\begin{array}{ccc} 1&3&\mathbf{5}\\2&\mathbf{6}\\\mathbf{4} \end{array}$$
where the latter pair is obtained from the former by applying a
single dual Knuth relation on the triple $\{4,5,6\}$.

$$\begin{array}{ccc} 1&2&\mathbf{4}\\\mathbf{3}&\mathbf{5}&6 \end{array}
\leq_{weak}
\begin{array}{ccc} 1&2&\mathbf{4}\\\mathbf{3}&6\\\mathbf{5} \end{array}
\hskip .1in \text{ but } \hskip .1in
\begin{array}{ccc} 1&2&\mathbf{3}\\\mathbf{4}&\mathbf{5}&6
\end{array} \not\leq_{weak}
\begin{array}{ccc} 1&2&\mathbf{5}\\\mathbf{3}&6\\\mathbf{4} \end{array}$$
where the latter pair is obtained from the former by applying a
single dual Knuth relation on the triple $\{3,4,5\}$.
\end{example}

\subsection{The definition of the skew orders}
\noindent

Let $m=k+n$, $\lambda \models m$ and $\mu\models k$ such that
$\mu\subset\lambda$. For  $T \in SYT_m$ of shape $\lambda$,  define
      $$ T_{\lambda/\mu}$$
to be  the skew standard tableau on $[n]$ of shape $\lambda/\mu$
obtained by standardizing the skew segment of $T$ having shape
$\lambda/\mu$.

\begin{definition}
\label{skew-orders-definition} Let $\leq$ be partial order on
$SYT_n$ having inner translation property. For  $U$ and  $V$ be two
skew standard tableaux  in $SYT_n^{\mu}$, we set
$$U \leq V $$
if there exist two tableaux $S$ and $T$ in $SYT_m$ of shape
$\lambda$ and $\lambda'$ respectively which  satisfy:
   $$\begin{aligned} & S_{\mu}=T_{\mu}=R \mbox{ for some } R \in SYT_k
                    \mbox{ of shape } \mu \\
          & S_{\lambda/\mu}=U \mbox{ and } T _{\lambda'/\mu}=V\\
          & S \leq T.
   \end{aligned}$$
\end{definition}

\begin{remark}\label{KL-geom-skew-poset-remark}
As a consequence of Theorem~\ref{Vogan-involution}, the skew orders,
$\leq_{KL}^{op}$ and $\leq_{geom}$ on $SYT_n^{\mu}$ becomes well
defined.
\end{remark}

\vskip .3in

\section{Proof of Theorem~\ref{skew-homotopy-theorem}}
\label{skew-homotopy-theorem-section}
 \vskip .1in

 In what follows we first  prove a result, namely
 Proposition~\ref{skew-mobius-value-theorem} below, which is
about  the M\"obius function of the subposet $SYT_m^R$~ ~of~
~$SYT_n$ ordered by any order that is stronger than  ~$\leq_{weak}$,
restricts to segments  and  has the property of extension from the
segments.  Consequently  Theorem~\ref{skew-homotopy-theorem} follows
from this results  together with
Theorem~\ref{embedding-from-initial-segments-theorem},
Theorem~\ref{Vogan-involution} and
Definition~\ref{skew-orders-definition}.

\vskip .1in
 Let  $R$ be a tableau in $SYT_k$ and $m=k+n$. Recall that
$$SYT_m^R:=\{T\ \in SYT_m \mid T_{[1,k]}=R\}.$$

Since the weak order restricts to segment, it  can be induced on
$SYT_m^R$. Moreover  the analysis made by comparing the left
inversion sets yields that any tableau $T \in SYT_m^R$, under the
weak order, lies between two tableaux $\hat{0}_{R,n}$ and
$\hat{1}_{R,n}$ given below.

 $\put(0,100){\makebox(1,1){}}
 \put(42,50){\makebox(12,35){$\hat{0}_{R,n}$=}}
 \put(85,65){$R$}
\put(65,90){\line(1,0){70}}
 \put(105,60){\line(1,0){30}}
\put(65,40){\line(1,0){40}} \put(135,60){\line(0,1){30}}
\put(105,40){\line(0,1){20}} \put(65,40){\line(0,1){50}}
\put(135.5,74.5){\framebox(57,15){{\scriptsize $k+1~\ldots~k+n$}}}
\put(232,50){\makebox(12,25){$\hat{1}_{R,n}$=}} \put(275,65){$R$}
\put(255,90){\line(1,0){70}}
 \put(295,60){\line(1,0){30}}
\put(255,40){\line(1,0){40}}
 \put(325,60){\line(0,1){30}}
\put(295,40){\line(0,1){20}}
 \put(255,40){\line(0,1){50}}
\put(255.5,4){\framebox(19,36){\shortstack{{\scriptsize
$k+1$}\\\\.\\.\\.\\{\scriptsize $k+n$}}}} $

\vskip .1in

 Therefore $ ( SYT_m^R, \leq_{weak}
)=[\hat{0}_{R,n},\hat{1}_{R,n}]_{\leq_{weak}}$  and for  any order
$\leq$ which is stronger than the weak order and which restricts to
segments we have
$$
(SYT_m^R, \leq)=[\hat{0}_{R,n},\hat{1}_{R,n}]_{\leq}
$$

\vskip .1in
\begin{proposition}
\label{skew-mobius-value-theorem} Let  $\leq$ be  any order on
$SYT_m$ with the following properties
\begin{enumerate}
\item[(i)] $\leq$ is stronger than $\leq_{weak}$
\item[(ii)] $\leq$ restricts to segments
\item[(iii)] $\leq$ extends from segments.
\end{enumerate}
 Then for  $\hat{0}_{R,n}$ and  $\hat{1}_{R,n}$  as above one has
$$ SYT_m^R=[\hat{0}_{R,n},\hat{1}_{R,n}],$$
 and
 the  proper part of $SYT_m^R$  is homotopy equivalent to
 $$
\begin{cases}
\text{ an } (n-2)-\text{dimensional sphere} & \text{ if } R \text{ is
rectangular,} \\
\text{ a point }                          & \text{ otherwise.} \\
\end{cases}
$$
\end{proposition}

\vskip .1in
Below we recall Rambau's  Suspension Lemma about bounded
posets \cite{Rambau}, which will be used  to prove
Proposition~\ref{skew-mobius-value-theorem}.

\begin{lemma}[Rambau's Suspension Lemma]
\label{suspension-lemma} Let $\mathcal{P}$ and  $\mathcal{Q}$ be two
bounded posets such that $\hat{0}_{\mathcal{Q}} \not =
\hat{1}_{\mathcal{Q}}$. Assume  $\mathcal{P}$ is the  disjoint union
of its two subsets $\mathcal{I}$ and $\mathcal{J}$ where
$\mathcal{I}$ forms an order ideal and $\mathcal{J}$ forms an order
filter of $\mathcal{P}$. Assume further that there are order
preserving maps
$$ f: \mathcal{P} \mapsto \mathcal{Q}~ \mbox {and } ~i,j: \mathcal{Q} \mapsto
\mathcal{P}$$ satisfying the following properties:
\begin{enumerate}
\item[(i)] The image of ~$i$~ lies in $\mathcal{I}$ and the image of ~$j$~ lies
in $\mathcal{J}$
\item[(ii)] The maps $f\circ i$ and $f\circ j$ are identity on $\mathcal{Q}$
\item[(iii)] For every $p \in \mathcal{P}$, $i\circ f(p) \leq p \leq j\circ
f(p)$
\item[(iv)] The fiber $f^{-1}(\hat{0}_{\mathcal{Q}})$ lies in $\mathcal{J}$
except for
$\hat{0}_{\mathcal{P}}$ and the fiber
$f^{-1}(\hat{1}_{\mathcal{Q}})$ lies in $\mathcal{I}$ except for
$\hat{1}_{\mathcal{P}}$.
\end{enumerate}

Then the proper part $\mathcal{P} -\{ \hat{0}_{\mathcal{P}},
\hat{1}_{\mathcal{P}}\}$ of $\mathcal{P}$ is homotopy equivalent to
the suspension of the proper part   of $\mathcal{Q}$.
\end{lemma}

\vskip .1in
 \begin{proof}[Proof of
 Proposition~\ref{skew-mobius-value-theorem}]
 For $n\geq2$, let
 $$
 \mathcal{P}=[\hat{0}_{R,n},\hat{1}_{R,n}]~ \mbox{ and }
  ~\mathcal{Q}=[\hat{0}_{R,n-1},\hat{1}_{R,n-1}]
  $$
together with the subposets of  $\mathcal{P}$ given as
$$
\begin{aligned}
\mathcal{I}=\{T \in \mathcal{P}: m-1 \not \in\Des(T)\}\\
\mathcal{J}=\{T \in \mathcal{P}: m-1  \in\Des(T)\}.
\end{aligned}
$$
Moreover let $$ f: \mathcal{P} \mapsto \mathcal{Q}~~\mbox{ and }~~
i,j: \mathcal{Q} \mapsto \mathcal{P} $$ where the map   $f$
restricts any $T \in \mathcal{P}$ to its initial segment
$T_{[1,m-1]}$ and the map $i$ concatenates $m$  to  the first row of
any $S\in \mathcal{Q}$ from right whereas  $j$ concatenates $m$ to
the first column of $S$ from the bottom.

 First  we will show
that $\mathcal{I}$ is an order ideal of $\mathcal{P}$.  Let  $T \in
\mathcal{I}$ and $T' < T$.  Then by
Lemma~\ref{tableau-descent-lemma}, $\Des(T') \subseteq \Des(T)$ and
 therefore $m-1$ does not  belong to $\Des(T')$. This shows that  $T' \in
\mathcal{I}$
 and $\mathcal{I}$ is an order ideal. A similar  argument also shows
 that $\mathcal{J}$ is an order filter of $\mathcal{P}$. On the other hand
 it can be easily seen that  $\mathcal{P}$ is the disjoint union of
$ \mathcal{I}$ and $\mathcal{J}$.

Since the tableau $R$ is common for both $\mathcal{P}$ and
$\mathcal{Q}$ and $\leq$  restricts to the initial segments, the map
$f: \mathcal{P} \mapsto \mathcal{Q}$ is well defined and order
preserving. By virtue of their definitions the maps $i,j:\mathcal{Q}
\mapsto \mathcal{P}$ are also well defined. On the other hand since
$\leq$ has the property of  extension from segments therefore they
both are order preserving.

Now part $(\mathrm{i})$ follows from the fact that the map $i$
concatenates $m$ to the right of the first row of $S \in
\mathcal{Q}$, which provides no possibility that $m$ appears below
$m-1$ in $i(S)$. Therefore $m-1 \not \in \Des(i(S))$ and $i(S) \in
\mathcal{I}$. On the other hand in $j(S)$, $m$ always appears below
$m-1$ and this shows that $j(S) \in \mathcal{J}$.

For part $(\mathrm{iii})$, let $\rho_T=a_{1}\ldots a_{l-1}~m~
a_{l+1} \ldots a_m$ be the row word of $T \in \mathcal{P}$. The
analysis on the
 (left) inversion sets gives:
$$a_{1}\ldots a_{l-1}~a_{l+1}\ldots a_m~m~\leq_{weak}~
 a_{1}\ldots a_{l-1}~m~a_{l+1} \ldots
a_m~ \leq_{weak}~ m~a_{1}\ldots a_{l-1}~a_{l+1} \ldots a_m$$ and by
$RSK$ correspondence $i\circ f(T) \leq_{weak} T \leq_{weak} j\circ
f(T)$ and hence $i\circ f(T) \leq T \leq j\circ f(T)$.

One can check the hypotheses $(\mathrm{ii})$ and $(\mathrm{iv})$
easily. Therefore by Lemma~\ref{suspension-lemma}, the proper part
of $\mathcal{P}$ is homotopy equivalent to the suspension of the
proper part of $\mathcal{Q}$.

\vskip .1in In the rest we proceed by induction: Let $n=1$. Then all
tableaux in the poset $\mathcal{P}=[\hat{0}_{R,1},\hat{1}_{R,1}]$
are obtained by placing $m$ in some outer corner of $R$, i.e, in an
empty cell along the boundary of $R$ whose addition to $R$ still
gives a Young tableau shape. Moreover it can be easily checked, for
example by comparing the left inversion sets of their row words,
that  these tableaux  form a saturated chain in $(SYT_m,
\leq_{weak})$. On the other hand since $\leq$ is stronger then the
$\leq_{weak}$ and  restricts to  segments this chain remains
saturated in $(SYT_m, \leq)$. The following diagram illustrates the
case when $R$ has three outer corners. \vskip .1in

$\put(0,80){\makebox(1,1){}}
 \put(22,30){\makebox(12,35){$\hat{0}_{R,n}$=}}
\put(150,30){\makebox(12,35){$\leq$}}
\put(290,30){\makebox(12,35){$\leq$}}
 \put(65,45){$R$}
\put(45,70){\line(1,0){70}}
 \put(85,40){\line(1,0){30}}
\put(45,20){\line(1,0){40}} \put(115,40){\line(0,1){30}}
\put(85,20){\line(0,1){20}} \put(45,20){\line(0,1){50}}
\put(115,60){\framebox(11,9.75){{\footnotesize $m$}}}
\put(215,45){$R$} \put(195,70){\line(1,0){70}}
 \put(235,40){\line(1,0){30}}
\put(195,20){\line(1,0){40}}
 \put(265,40){\line(0,1){30}}
\put(235,20){\line(0,1){20}}
 \put(195,20){\line(0,1){50}}
\put(235,30){\framebox(11,9.75){{\footnotesize $m$}}}
\put(418,33){\makebox(12,25){=$\hat{1}_{R,n}$}} \put(355,45){$R$}
\put(335,70){\line(1,0){70}}
 \put(375,40){\line(1,0){30}}
\put(335,20){\line(1,0){40}}
 \put(405,40){\line(0,1){30}}
\put(375,20){\line(0,1){20}}
 \put(335,20){\line(0,1){50}}
\put(335,10){\framebox(11,9.75){{\footnotesize $m$}}} $

Now if $R$ has rectangular shape then it has two outer corners and
the poset $\mathcal{P}=[\hat{0}_{R,1},\hat{1}_{R,1}]$ consists of
two tableaux. It has  the M\"obius function from the bottom to the
top elements to be $-1$ and moreover the proper part of
$\mathcal{P}$ is homotopy equivalent to  the empty set i.e,
$(-1)$-dimensional sphere.

If $R$ is non rectangular then as in the above diagram $\mathcal{P}$
is a saturated chain having more than two elements. Hence its
M\"obius function is $0$ from  the bottom to the top elements and it
is homotopic to  a point.

Now assume that for $n=r$ the poset
$\mathcal{Q}=[\hat{0}_{R,r},\hat{1}_{R,r}]$ satisfies the hypothesis
i.e., the proper part of $\mathcal{Q}$ is homotopic to a
$(r-2)$-sphere in case $R$ is rectangular and it is homotopic to a
point otherwise.

On the other hand we already see  that the proper part
 of $\mathcal{P}=[\hat{0}_{R,r+1},\hat{1}_{R,r+1}]$ is homotopy
equivalent to the suspension of the proper part
 of $\mathcal{Q}$, so that the
former becomes homotopy equivalent to a ($r-1$)-sphere if $R$ is
rectangular and to a point otherwise. Therefore the assertion of
Proposition~\ref{skew-mobius-value-theorem} follows.
\end{proof}

\vskip .1in

\begin{proof}[Proof of Theorem~\ref{skew-homotopy-theorem}] By
Remark~\ref{KL-geom-skew-poset-remark}, $KL$ and geometric orders
are well defined on $SYT_n^{\mu}$. On the other hand they restrict
segments and  have the property of embedding from initial segments
by Lemma~\ref{embedding-from-initial-segments-theorem}. So the
required statement follows from
Proposition~\ref{skew-mobius-value-theorem}.
\end{proof}

\vskip .3in
\section{Shortest and longest chains}
\label{Shortes-Longest-Chains-section}
\vskip .1in
 By observing Figure \ref{figure1}, one can see
that the posets of $SYT_n$ with all these orders are not lattices
and not ranked. On the other hand we can still say something about
the size of their shortest or longest chains, where by convention
$c_1 < c_2 < \ldots < c_i$ has size $i$.

\begin{proposition}
\begin{enumerate}
\item[]
\item[(i)] The size of a  shortest saturated chain in $(SYT_n,
\leq_{weak})$ is $n$.
\item[(ii)] The size of a longest chain  in
$(SYT_n, \leq_{weak})$, $(SYT_n, \leq^{op}_{KL})$ and  $(SYT_n,
\leq_{geom})$ is equal to the size of the longest chain in $(\Par_n,
\leq_{dom})$, which is asymptotically $(\sqrt{8}n^{3/2})/{3}$.
\end{enumerate}
\end{proposition}

\begin{proof}
Observe that if $\sigma$ is covered by $\tau$ in $(\sym_n,
\leq_{weak})$ then the size of the (left) descent set $\Des_L(\tau)$
of $\tau$ is at most one bigger than the size of $\Des_L(\sigma)$.
This fact is also true for $(SYT_n, \leq_{weak})$: If $S$ is covered
by $T$ in $(SYT_n, \leq_{weak})$  then
\begin{equation} \label{covering-descent-relation2}
0~\leq~\mid\Des(T)~\backslash~\Des(S)\mid~<~2~.
\end{equation}
This shows that the size of a shortest saturated chain must be  at
least $n$. On the other hand it can be seen  by  an easy induction
that there exist a saturated chain in $(SYT_n, \leq_{weak})$ of size
$n$ with the following form:
\begin{equation}\label{a-shortest-chain}
\begin{array}{ccccc} 1&2&3\ldots&n\end{array}~\leq~
\begin{array}{ccccc} 1&2&3\ldots&n-1\\n&&&&\end{array}\leq~
~\ldots~ \leq~
\begin{array}{cc} 1&2\\3&\\ \vdots &\\ n& \end{array}~\leq~
\begin{array}{c} 1\\ 2\\3\\\vdots\\n \end{array}.
\end{equation}
Therefore the statement about shortest chains in $(SYT_n,
\leq_{weak})$ follows.

 For longest chains the  proof is based on two facts: a result
 of  Greene and Kleitman \cite[Page~9]{Greene&Kleitman} which calculates the size of longest
chain in the lattice of integer partitions ordered by the dominance
order and the result of  Melnikov  \cite[Proposition
4.1.8]{Melnikov1} which shows that for any tableau $S$  of shape
$\mu$  in $SYT_{n}$ and  for any partition $\lambda \models n$ such
that  $\mu <_{dom}^{op} \lambda$, there is a tableau $T \in SYT_n$
such that $\sh(T)=\lambda$  and  $S <_{weak}T$. These two facts
enable us to calculate the longest chain in $SYT_{n}$ ordered by the
weak order. Since $ \leq^{op}_{KL}$ and $\leq_{geom}$ also change
the shapes of the tableaux, the longest chain of $(SYT_n,
\leq_{weak})$ still remains saturated in $(SYT_n, \leq^{op}_{KL})$
and $(SYT_n, \leq_{geom})$.
\end{proof}

\begin{remark}
By an easy induction one can see that chain in
(\ref{a-shortest-chain}) still remains  saturated in  $SYT_n$ for
$KL$, geometric and chain orders. Therefore if it were known
(\ref{covering-descent-relation2}) is satisfied by these three
orders, we could deduce the same conclusion about their shortest
chains.
\end{remark}

\vskip .2in
\section{Remarks and questions}
\label{remarks-section}
\vskip .1in

\begin{remark}
Theorem~\ref{homotopy-theorem} also follows from Proposition
~\ref{skew-mobius-value-theorem} by taking $R= 1$. The original
proof is kept here  for indicating different approaches to the
subject.
\end{remark}

\begin{remark}
The order complex of the proper part of  $(SYT_n,\leq)$ under any of
the four orders is not homeomorphic to a sphere. One can observe
$SYT_4$ in Figure \ref{figure1} to see the smallest example.
Moreover since these posets are not ranked for $n\geq4$, the order
complex of their proper parts are not pseudomanifolds.
\end{remark}

\begin{remark}
\label{inner-translation-remark}
Although the weak order on $SYT_n$ does not have the inner
translation property, it might still satisfy
Corollary~\ref{equal-inner-shapes} without this property, which
would then make it possible to define weak  order on skew standard
tableaux.

 For chain order, two pairs of tableaux given below where the inner
tableau {\small $\begin{array}{cc}1&2\\3\end{array}$} common  to the
first pair is  replaced by
$\scriptstyle{\begin{array}{cc}1&3\\2\end{array}}$ in the second
pair, show that Corollary~\ref{equal-inner-shapes} is not satisfied
by chain order:
$$
\begin{array}{ccc} \textbf{1}&\textbf{2}&5\\\textbf{3}&6\\4 \end{array}
\geq_{chain}
\begin{array}{ccc} \textbf{1}&\textbf{2}&5\\\textbf{3}&4\\6 \end{array} \hskip
.2in
\text{but} \hskip .2in
\begin{array}{ccc} \textbf{1}&\textbf{3}&5\\\textbf{2}&6\\4 \end{array} \not
\geq_{chain}
\begin{array}{ccc} \textbf{1}&\textbf{3}&5\\\textbf{2}&4\\6 \end{array}
$$
\end{remark}

\begin{question}
One might ask to what extent the definitions and results in this
paper apply to other Lexicographic Coxeter systems $(W,S)$.  The
weak order on $W$ is well-defined, as are $KL$
 and the
geometric order, where the former still remains weaker than the
latter (\cite{Kazhdan-Lusztig}; see \cite[Fact 7]{Garsia-McLarnan}).
Definition~\ref{weak-order-def} makes sense and remains valid, and
so does Proposition~\ref{order-preserving-maps}$\mathrm{(i)}$ for
$KL$ order (\cite{Kazhdan-Lusztig}; see \cite[Fact
7]{Garsia-McLarnan}). For geometric order  the same property follows
from \cite[Theorem 6.11]{Borho} or \cite[Theorem 9.9]{Joseph2}.

 For the analysis of M\"obius function and homotopy
types, the crucial Lemma~\ref{descent-fibers} was proven by Bjorner
and Wachs \cite[Theorem 6.1]{Bjorner3}  for all {\it finite} Coxeter
groups $W$. Hence Corollary~\ref{SYT-homotopy-corollary} and
Theorem~\ref{homotopy-theorem} are valid also in this generality,
with the same proof.
\end{question}

\section*{Acknowledgements}
The author is grateful  to  V. Reiner for his helpful questions and
comments throughout this work. The author also thanks A. Melnikov
for allowing access to her unpublished preprints, W. McGovern and M.
Geck  for helpful comments.

\end{document}